\documentclass{article}
\usepackage{amsfonts}
\usepackage{amssymb}
\usepackage{graphicx}
\usepackage{amsmath}

\newtheorem{theorem}{Theorem}

\newtheorem{conjecture}[theorem]{Conjecture}
\newtheorem{corollary}[theorem]{Corollary}

\newtheorem{definition}[theorem]{Definition}

\newtheorem{proposition}[theorem]{Proposition}
\newtheorem{remark}[theorem]{Remark}

\newenvironment{proof}[1][Proof]{\textbf{#1.} }{\  \rule{0.5em}{0.5em}}
\begin{document}

\title{A new proof of the virtual Haken conjecture}
\author{Charalampos Charitos}
\maketitle

\begin{abstract}
A new direct proof of the Virtual Haken Conjecture, which asserts that every
compact, orientable, irreducible three-dimensional manifold with infinite
fundamental group has a finite cover that is Haken, will be given.

\medskip

\bigskip

\textbf{Keywords}, Incompressible surfaces, Haken 3-manifolds
\end{abstract}

\section{Introduction}

A closed, orientable 3-manifold $M$ is said to be Haken if it contains an
incompressible embedded surface. The virtually Haken conjecture asserts that
every closed, orientable, irreducible 3-manifold with infinite fundamental
group is virtually Haken; that is, it possesses a finite cover that is
Haken. The conjecture is commonly attributed to Waldhausen \cite{Waldhausen}%
, although he did not state it explicitly and it appears as Problem 3.2 in
Kirby's problem list \cite{Kirby}. Thurston emphasized the virtually Haken
conjecture as a cornerstone of 3-manifold theory. He included it as the
sixteenth problem in his well-known list of twenty-four problems \cite%
{Thurston1}, where he outlined his program for understanding 3-manifolds.
Following Perelman's proof of Thurston's geometrization conjecture, the
virtually Haken conjecture remained unresolved only for hyperbolic
3-manifolds. Agol subsequently established the conjecture, building on
methods developed by Wise. Perelman's theorem thus represents a crucial link
in the chain of results leading to the proof of the virtually Haken
conjecture. Nevertheless, since Perelman's proof relies on analytic
techniques and Riemannian geometry, a direct proof of the virtually Haken
conjecture using purely geometric and topological arguments, in the spirit
of Thurston's ideas, remains an important challenge.

This paper is organized as follows. Section 2 reviews Thurston's
hyperbolization theorems and introduces the notion of atoroidal 3-manifolds.

Section 3 states a theorem of Cooper, Long, and Reid \cite{Cooper-Long-Reid}%
, which is central to our approach, as it motivates the strategy adopted
here. In particular, elements of the proof of Theorem 1.1 in \cite%
{Cooper-Long-Reid} are adapted to establish our main result.

Section 4 proves that every closed, irreducible, atoroidal 3-manifold with
infinite fundamental group can be obtained by gluing a handlebody $H$ of
genus $g\geq 2$ to a manifold $N$ that is irreducible, boundary-irreducible,
atoroidal, and acylindrical, via a homeomorphism $f:\partial N\rightarrow
\partial H.$

Finally, Section 5 contains the proof of the virtually Haken conjecture.

\section{Thurston's great theorems}

This section states Thurston's hyperbolization theorems, which will be used
in the sequel and which have been fundamental to the development of
3-manifold theory.

\subsection{Hyperbolic manifolds}

\begin{theorem}
\label{Hyperbolization for closed} Let $M$ be a closed orientable manifold.
Then $M$ admits a hyperbolic structure provided that:

$(1)$ $M$ is irreducible and Haken;

$(2)$ $M$ is atoroidal i.e. $\pi _{1}(M)$ does not contain a subgroup
isomorphic to $\mathbb{Z}\oplus \mathbb{Z}.$
\end{theorem}

Assuming that $M$ admits a hyperbolic structure it follows that $M=H^{3}/G,$
where $G$ is a discrete group of isometries of the hyperbolic space $H^{3}$
which acts cocompactly on $H^{3}.$ Therefore, $M$ must be irreducible. Also,
since $M$ is compact, there is an $\varepsilon >0$ such that any closed loop
of length $\leq \varepsilon $ is homotopically trivial. Therefore, no
nontrivial element of $G=\pi _{1}(M)$ can move any point of $H^{3}$ at a
distance less than $\varepsilon .$ Parabolic isometries do not have this
property, so all elements of $G$ are hyperbolic isometries. Now, let $\Gamma 
$ be an abelian subgroup of $G$ and let $g\in \Gamma .$ Then $g$ fixes an
axis $L_{g}$ and since $\Gamma $ is abelian any other element of $\Gamma $
must fix the same axis. Thus $\Gamma $ must be cyclic of infinite order.
These observations show that it is essential to assume that $M$ is both
irreducible and atoroidal. By contrast, there exist many hyperbolic,
non-Haken 3-manifolds. William Thurston demonstrated that all but finitely
many Dehn fillings on the figure-eight knot complement produce irreducible,
non-Haken, non-Seifert-fibered manifolds---and moreover, these manifolds are
hyperbolic. This striking result provided the first infinite family of such
examples. Prior to Thurston's work, non-Haken manifolds were thought to be
exceedingly rare: the only known examples were certain Seifert fibered
spaces over $S^{2}$ with three exceptional fibers.

\begin{theorem}
\label{Hyperbolization for cusps}Assume that $M$ is a compact, irreducible,
orientable $3$-manifold with $\partial M\neq \emptyset .$ Assume that $%
\partial M$ consists of tori and that any abelian non-cyclic subgroup $%
\Gamma $ of $\pi _{1}(M)$ is peripheral i.e. $\Gamma $ is conjugate to the
fundamental subgroup of a component of $\partial M.$ Then $Int(M)$ admits a
hyperbolic structure of finite volume, provided that $M$ is not $S^{1}\times
D^{2},$ $T^{2}\times I$ and the orientable $I$-bundle over the Klein bottle.
\end{theorem}

Concerning Theorem, Thurston constructed many beautiful examples of
finite-volume hyperbolic structures on the complements of various knots and
links in $S^{3}.$ Moreover, he proved that every knot and link in $S^{3}$ is
exactly one of three mutually exclusive types---hyperbolic, torus, or
satellite---and that precisely those that are neither torus nor satellite
have hyperbolic complements. These results provided a vast, explicit family
of hyperbolic complements, \cite{Thurston} or (see \cite{Marden} p. 325).

Assuming now that $\partial M\neq \emptyset $ and each component of $%
\partial M$ is of genus $\geq 2$ we have the following theorem:

\begin{theorem}
\label{Hyperbolization for compact} Let $M$ be an orientable, compact,
irreducible, atoroidal 3-manifold. We assume that $\partial M\neq \emptyset $
and that $M$ is boundary irreducible i.e. every component of $\partial M$ is
incompressible and $M$ does not contain an essential properly embedded
annulus. Then $M$ admits a hyperbolic structure with geodesic boundary.
\end{theorem}

In connection with Theorem \ref{Hyperbolization for compact}, Thurston
constructed striking examples of finite-volume hyperbolic structures on
compact 3-manifolds with geodesic boundary by gluing together truncated
hyperideal hyperbolic polyhedra (see Example 3.2.12, p. 133 in \cite%
{Thurston 2}). Moreover, if $M$ admits a hyperbolic structure with geodesic
boundary, then $M$ cannot contain any essential annulus. Indeed, doubling $M$
along its geodesic boundary produces a closed manifold $DM$ that is
hyperbolic and thus cannot contain an essential torus which is forbidden by
the remarks following Theorem \ref{Hyperbolization for closed}.

Finally, for the interior of compact manifold the following theorem is valid
(see Theorem 2.3 in \cite{Thurston1}).

\begin{theorem}
\noindent \label{Hyperbolization for the interior} The interior of a compact
3-manifold $M$ with nonempty boundary has a hyperbolic structure if and only
if $M$ is irreducible, atoroidal and not homeomorphic to the orientable $I-$%
bundle over the Klein bottle.
\end{theorem}

Generally, in the case where $\partial M\neq \emptyset ,$ it is well-known
that $M$ contains properly embedded incompressible surfaces. Indeed, the
presence of such surfaces is indispensable for the proofs of the theorems
discussed above. By a theorem of Waldhausen, elaborated for example in
Hempel \cite{Hempel}, $M$ admits a hierarchy, meaning it can be successively
cut along incompressible surfaces into simpler 3-manifolds until one ends up
with a disjoint union of 3-balls. This hierarchical decomposition allows one
to construct hyperbolic structures inductively on the simplest pieces and
then glue them together systematically to obtain a hyperbolic structure on
the original manifold, following Thurston's approach.

We conclude this discussion by stating Thurston's hyperbolization
conjecture, which was proven by Grigori Perelman in 2003, building upon
Richard Hamilton's work on the Ricci flow.

\begin{conjecture}
Each closed, orientable, atoroidal 3-manifold with infinite fundamental
group admits a hyperbolic structure.
\end{conjecture}

It is worth noting here that the previous conjecture can be proven assuming
that the virtual Haken conjecture holds.

\subsection{The atoroidality condition}

In this paragraph we will clarify the term \textquotedblleft atoroidal
manifold\textquotedblright , as its definition varies across different
contexts. In 3-manifold topology, a manifold is termed atoroidal if it does
not contain any essential tori. An essential torus is typically defined as
an embedded, incompressible torus that is not boundary-parallel. However,
the precise definition can differ: some authors adopt a geometric
perspective, while others use an algebraic approach based on the fundamental
group. This variation in definitions can lead to different interpretations
of what constitutes an atoroidal manifold. Therefore, it is crucial to
specify the intended definition when discussing atoroidal manifolds to avoid
ambiguity.

\begin{definition}
(i) A closed 3-manifold $M$ is called \emph{homotopically atoroidal} if
there is no map $T\rightarrow M$ from the 2-torus $T$ to $M$ which induces a
monomorphism $\pi _{1}(T)\rightarrow \pi _{1}(M).$ If $M$ is compact with $%
\partial M\neq \emptyset $ then $M$ is called \emph{homotopically atoroidal}
if the existence of a map $T\rightarrow M$ which induces a monomorphism $\pi
_{1}(T)\rightarrow \pi _{1}(M)$ implies that the map can be homotoped into a
boundary component of $M.$

This condition can be rephrased easily in terms of the fundamental group.
That is, a 3-manifold $M$ is called homotopically atoroidal if each rank 2
abelian subgroup of $\pi _{1}(M)$ is conjugate to $\pi _{1}(\partial _{0}M),$
where $\partial _{0}M$ is some toral boundary component of $\partial M.$

(ii) A closed 3-manifold $M$ is called \emph{geometrically atoroidal} if
there is no embedding $T\rightarrow M$ which induces a monomorphism $\pi
_{1}(T)\rightarrow \pi _{1}(M).$ If $M$ is compact with $\partial M\neq
\emptyset $ then $M$ is called \emph{geometrically atoroidal} if the
existence of an embedding $T\rightarrow M$ which induces a monomorphism $\pi
_{1}(T)\rightarrow \pi _{1}(M)$ implies that the image of $T$ can be
isotoped into the boundary of $M$ i.e. the image of $T$ is parallel to a
boundary component of $\partial M.$
\end{definition}

Recall that a small Seifert manifold is a Seifert manifold which does not
contain a fiberwise incompressible torus. In \cite{Jaco} small Seifert
manifolds are referred to as special manifolds and they are described in
page 155 of \cite{Jaco}. The following result is known as the Torus Theorem,
and it asserts that:

\begin{theorem}
\label{torus theorem} Let $M$ be a compact orientable irreducible
3-manifold. We also assume that $M$ is not a small Seifert manifold. Then $M$
is homotopically atoroidal if and only if $M$ is geometrically atoroidal.
\end{theorem}

The proof of the Theorem \ref{torus theorem} is the culmination of the work
of many mathematicians (see for instance Theorem 1.40 and 1.41 in \cite%
{Kapovich} and the Paragraph 1.6 in \cite{AFW}, as well as, the citations
and the scholia accompanying these theorems). More specifically, if $%
\partial M\neq \emptyset ,$ Theorem \ref{torus theorem} is the torus
theorem, see \cite{Scott1} and Theorem VIII.14 in \cite{Jaco}, while if $%
\partial M=\emptyset $ the JSJ (Jaco, Shalen, Johannson) decomposition
theorem \cite{Jaco-Shalen}, \cite{Johansson}, \cite{Neumann} is used to
accomplish the proof \cite{Scott2}.\smallskip

\noindent \textbf{Convention}. In what follows by the term \textit{atoroidal
manifold} we will mean that the manifold is geometrically atoroidal.

\section{The virtual Haken conjecture for 3-manifolds with boundary}

We begin this section by presenting a theorem by D. Cooper, D. Long, and A.
Reid \cite{Cooper-Long-Reid}, which provided a significant impetus toward
proving the Virtual Haken Conjecture. This theorem also offers a resolution
to the conjecture in the context of 3-manifolds with boundary.

\begin{definition}
Let $M$ be a compact orientable 3-manifold. A map $S\rightarrow M$ of a
closed orientable connected surface surface $S$ is \emph{essential} if it is
injective at the level of fundamental group and the group $i_{\ast }\pi
_{1}(S)$ cannot be conjugate to a subgroup $\pi _{1}(\partial _{0}M)$ of $%
\pi _{1}(M)$, where $\partial _{0}M$ is a component of $\partial M.$

By abusing the language the surface $S$ will be also called \emph{essential}
in $M.$
\end{definition}

Theorem 1.1 in \cite{Cooper-Long-Reid} asserts that:

\begin{theorem}
\label{essential1}Let $M$ be a compact connected 3-manifold with non-empty
incompressible boundary. Suppose that the interior of $M$ has a complete
hyperbolic structure. Then either $M$ is covered by $F\times I,$ where $F$
is a closed orientable surface, or $M$ contains an immersed essential
surface $S$ of genus at least $2.$

Furthermore $S$ can be lifted to an embedded nonseparating surface in a
finite cover of $M.$
\end{theorem}

The essential surface $S$ in Theorem \ref{essential1} is immersed in $M$
even if this in not mentioned in Theorem 1.1 of \cite{Cooper-Long-Reid}.
This follows from the fact that $S$ is the projection to $M$ of an embedded
surface $\widetilde{S}$ constructed in a finite cover of $M.$

From Theorem \ref{Hyperbolization for the interior} we get the following
result, as a corollary of the previous theorem.

\begin{theorem}
\label{essential2} Let $M$ be a compact, connected, orientable irreducible
and atoroidal 3-manifold with non-empty incompressible theorem. Then either $%
M$ is covered by a product $F\times I,$ where $F$ is a closed orientable
surface, or $M$ contains an essential closed surface $S.$

Furthermore $S$ can be lifted to an embedded non-separating surface in a
finite cover of $M.$
\end{theorem}

\begin{proof}
The manifold $M$ in Theorem \ref{essential2} admits a complete hyperbolic
structure in its interior thanks to Thurston's Theorem \ref{Hyperbolization
for the interior}. Thus, Theorem \ref{essential1} can be applied in order to
get the result.
\end{proof}

Some ideas underlying the proof of Theorem 1.1 \cite{Cooper-Long-Reid} will
be adapted in the final section of this work to establish the Virtual Haken
Conjecture.

\section{The complement of links in hyperbolic 3-manifolds}

Our first result in this paragraph is mentioned without proof in \cite%
{Marden}, p. 326, (see the remarks after the theorem referred to as
\textquotedblleft Dehn Surgeries on Hyperbolic Link
Complements\textquotedblright ).

Let $M$ be a closed orientable 3-manifold $M.$ It is a well known theorem
that $M$ can be obtained by surgery from a link $\mathcal{L}=L_{1}\cup
L_{2}\cup $ $\cdots \cup $ $L_{n}$ in $S^{3}$ which has the following
properties: \smallskip

\textit{(i) Each component }$L_{i}$\textit{\ is unknotted.}

\textit{(ii) Considering an orientation on }$L,$\textit{\ the linking number 
}$lk(L_{i},L_{j})$\textit{\ of any two components of }$L$\textit{\ is }$\pm
1 $\textit{\ or }$0.$

\textit{(iii) The surgery coefficient for each component }$L_{i}$\textit{\
is }$r_{i}=\pm 1.\smallskip $

A proof of this fundamental statement can be found in many sources (see for
instance \cite{Rolfsen}, Theorem 1, Chapter 9, p. 273 and Remark 8, p. 279).

Also for the link $\mathcal{L}$ we may assume that\smallskip

\textit{(iv) The number of components }$L_{i}$\textit{\ of }$\mathcal{L}$ 
\textit{such that }$\mathcal{L}$\textit{\ has properties (i)-(iii) is
minimum.\smallskip }

Below we assume that $M$ is additionally irreducible, atoroidal with
infinite fundamental group. Then, as a corollary to the following theorem,
we will demonstrate that $\mathcal{L}$ possesses the additional property:

\textit{(v) For each family of components }$\{L_{i_{1}},...,L_{i_{k}}\}%
\subseteq \{L_{1},...,L_{n}\}$\textit{\ there exist }$L_{i_{j}}\in
\{L_{i_{1}},...,L_{i_{k}}\}$\textit{\ and }$L_{j}\in \{L_{1},...,L_{n}\}$%
\textit{\ such that }$lk(L_{j},L_{i_{j}})\neq 0$\textit{\ and }$%
lk(L_{j},L_{i_{m}})=0$\textit{\ for each }$L_{i_{m}}\in
\{L_{i_{1}},...,L_{i_{k}}\} \backslash \{L_{i_{j}}\}.$

\begin{center}
\begin{figure}[h]
%\hspace*{9mm}
\includegraphics[bb = 43 486 513 699]{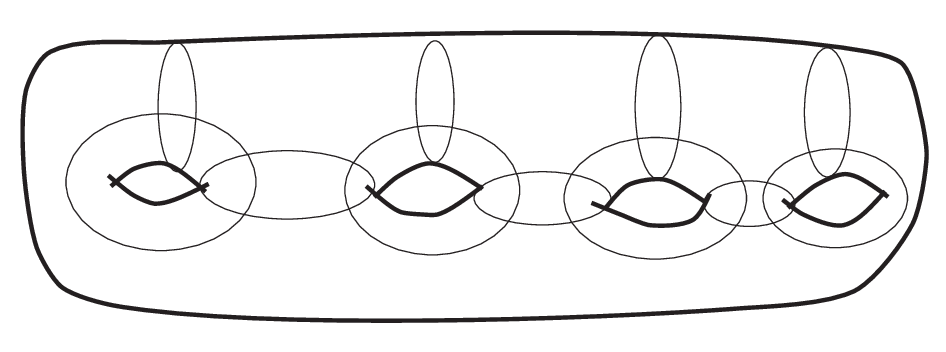} 
%\includegraphics[scale=0.74]{HAKEN_1.pdf} %\\vspace*{-8mm}\newline
%\label{disk4}
\caption{The curves $d_{i}.$}
\end{figure}
\end{center}

\begin{theorem}
\label{link complement}Every closed, irreducible, atoroidal 3-manifold $M$
with infinite fundamental group is obtained by Dehn surgery along some link $%
\mathcal{L}\subset S^{3}$ whose complement $N_{0}$ admits a hyperbolic
structure of finite volume. Furthermore, $N_{0}$ is acylindrical and each
boundary component of $N_{0}$ is incompressible.
\end{theorem}

\begin{proof}
We noted above that $M$ can be obtained by surgery from a link $\mathcal{L}%
=L_{1}\cup L_{2}\cup $ $\cdots \cup $ $L_{n}$ in $S^{3}$ which has
Properties (i)-(iv).

In what follows we will say that $M$ is obtained from $\mathcal{L}$ by
performing surgery with surgery coefficients $\mathbf{r}=(r_{1},\cdots
,r_{n}).$

From conditions (i) and (ii) we deduce that $\mathcal{L}$ must consist of at
least two components. Indeed, if $\mathcal{L}$ had only one component, then
the result of the surgery, $M,$ would be $S^{3}.$ In fact, one can show that 
$\mathcal{L}$ must consist of at least three components; however, this
conclusion is not immediate and is not used in the subsequent arguments.
Also, since $M$ is irreducible and hence prime, and in combination with
condition (iv), we deduce that there is not a component $L_{i_{0}}$ of $%
\mathcal{L}$ which has linking number $0$ with all the other components of $%
\mathcal{L}.$

For each $i,$ let $V_{i}$ be a solid torus around $L_{i};$ let $%
N_{0}=M\backslash (\cup _{i}Int(V_{i}))=S^{3}\backslash (\cup
_{i}Int(V_{i})).$ We will prove that $N_{0}$ is irreducible and atoroidal
and thus from Theorem \ref{Hyperbolization for cusps}, $Int(N_{0})$ admits a
hyperbolic structure. We will also prove that $N_{0}$ is acylindrical which
is another crucial property of $N_{0}.$ Finally we will show that each
boundary component of $\partial N_{0}$\ is incompressible in $N_{0}.$

In fact, we have:

(1) $N_{0}$\textbf{\ }\textit{is irreducible}.

Let $S$ be a 2-sphere embedded in $N_{0}.$ First we remark that $S$
separates $N_{0}.$ Indeed, if that were not the case, then $S$ would fail to
separate $S^{3}$, which is impossible. Let $C^{+}$ and $C^{-}$ be the two
components of $N_{0}\backslash S.$ Then all solid tori $V_{i}$ must be
contained in one of the components $C^{+}$ or $C^{-}.$ Indeed, if that were
not the case, then some components $L_{i}$ must belong in $C^{+}$ and some
others in $C^{-}.$ But then, either $M$ is not irreducible or we get a
contradiction in condition (iv). So, we may assume that all $V_{i}$ are
contained in $C^{+}.$ Therefore $C^{-}$ is homeomorphic to a 3-ball since $S$
lives in $S^{3}$ and thus $C^{-}$ coincides with one component of $%
S^{3}\backslash S.$ Therefore, $N_{0}$ is irreducible.

(2) $N_{0}$\textbf{\ }\textit{is atoroidal}.

Let $T$ be a torus in $N_{0}.$ First we will show that $T$ separates $N_{0}.$
Indeed, if $T$ does not separate $N_{0}$ then $T$ would not separate $S^{3}.$
But this is impossible by the solid torus theorem which asserts that each
torus in $S^{3}$ separates off a solid torus from at least one side (see for
example \cite{Rolfsen} Ch. 4, p. 107).

Let $R$ and $R^{\prime }$ denote the closures of the two components of $%
N_{0}\backslash T.$ First we assume that some boundary components of $N_{0}$
are contained in $R$ and some others in $R^{\prime }.$ Considering $T$ in $%
M, $ we remark that $T$ is still separating in $M$ and let us denote by $%
R_{M}$ and $R_{M}^{\prime }$ the closures of the two components of $%
M\backslash T$ corresponding to $R$ and $R^{\prime }.$ Obviously, some
components of $\mathcal{L}$ (which correspond to boundary components of $%
N_{0})$ are contained in $R_{M}$ and some others in $R_{M}^{\prime }.$ Since 
$M$ is atoroidal we deduce, without loss of generality, that $R_{M}^{\prime
} $ is a solid torus and let $V_{j}\subset R_{M}^{\prime },$ $j=1,\cdots ,k$
with $k\geq 1.$ The torus $T$ is also lying in $S^{3}$ and let us denote by $%
Q,$ $Q^{\prime }$ the closures of the components of $S^{3}\backslash T,$
where $L_{j}\subset Q^{\prime },$ $j=1,\cdots ,k.$ Let $\mathcal{L}^{\prime
}=L_{1}\cup \cdots L_{k}.$ By the solid torus theorem mentioned above at
least one of $Q,$ $Q^{\prime }$ is a solid torus. If $Q^{\prime }$ was a
solid torus then $M$ would be obtained by performing surgery along the link $%
L_{k+1}\cup \cdots \cup L_{n}$ with surgery coefficients $(r_{k+1},\cdots
,r_{n}).$ But this is impossible by the hypothesis (iv) above. So $Q$ must
be a solid torus and at the same time $Q^{\prime }$ is not a solid torus,
i.e. $Q^{\prime }$ is the complement of a non-trivial knot of $S^{3}.$ But
in that case we also reach a contradiction, otherwise, by performing surgery
along along $\mathcal{L}^{\prime }\subset Q^{\prime }$ we would obtain the
solid torus $R_{M}^{\prime },$ which is unacceptable.

From the preceding analysis, we may assume that all $V_{i}$ are contained in 
$R\subset N_{0}.$ Now $T,$ as a torus in the atoroidal manifold $M,$ bounds
a solid torus $B\subset M.$ If $B=R_{M}^{\prime }$ then $R^{\prime }$ is a
solid torus and hence $N_{0}$ is atoroidal. If $B=R_{M},$ we encounter a
contradiction as follows: gluing properly $B$ with a solid torus, say $%
B_{0}, $ we may take the sphere $S^{3}$ in which the link $\mathcal{L}$ is
contained and in which we have performed surgery with coefficients $\mathbf{r%
}=(r_{1},\cdots ,r_{n}).$ However, this leads to a contradiction, because
performing surgery along $\mathcal{L}$ with surgery coefficients $r$ should
yield the manifold $M.$

(3) $N_{0}$\textbf{\ }\textit{is acylindrical}.

In fact, let $Q$ be an annulus in $N_{0}$ with $\partial Q\subset \partial
N. $ Then there exists a torus $T$ in $N_{0}$ which is built either by $Q$
and one annulus $A\subset \partial N_{0}$ if $\partial Q$ is contained in
one component of $\partial N_{0}$ or, by two parallel copies $Q_{1},$ $Q_{2}$
of $Q$ and two annuli $A_{1},$ $A_{2}$ contained in two distinct components $%
T_{1},$ $T_{2}$ of $\partial N_{0}.$ Since $N_{0}$ is atoroidal $T$ must
bound a solid torus in $N_{0}.$ But in the former case, we deduce that $Q$
is parallel to a sub-annulus $A$ of $\partial N_{0}.$ In the latter case, we
deduce that $T_{1},$ $T_{2}$ must be parallel. Consequently, $M$ should be a
lens space. since $M$ is obtained by gluing two solid tori along their
boundaries. Therefore, we arrive at a contradiction, which completes the
proof of our statement.

(4) \textit{Each boundary component }$S_{i}$\textit{\ of }$\partial N_{0}$%
\textit{\ is incompressible}.

Let us assume for example that the component $S_{1}$ is compressible in $%
N_{0}.$ Then $S_{1}$ would bound a solid torus $W$ in $M$ containing all the
components $S_{2},\cdots ,S_{n}.$ Then we may find in $W$ a 3-ball $B$
containing all the components $L_{2},\cdots ,L_{n}$ of $\mathcal{L}.$
However, this is impossible because some of the components $L_{i},$ $i\neq 1$
must have a linking number of $\pm 1$ with $L_{1}.$
\end{proof}

\begin{corollary}
\label{hyperbolic link complement}Every closed, hyperbolic 3-manifold $M$ is
obtained by Dehn surgery along some link $\mathcal{L}\subset S^{3}$ whose
complement is hyperbolic.
\end{corollary}

\begin{proof}
Since $M$ is hyperbolic it follows that $M$ is atoroidal irreducible with
infinite fundamental group. Therefore by applying the the previous Theorem %
\ref{link complement} and Theorem \ref{Hyperbolization for cusps} we deduce
that $N_{0}=M\backslash (\cup _{i}Int(V_{i}))=S^{3}\backslash (\cup
_{i}Int(V_{i}))$ admits a complete finite-volume hyperbolic structure.
\end{proof}

\begin{corollary}
\label{property v} For each family of components $\{L_{i_{1}},...,L_{i_{k}}%
\} \subseteq \{L_{1},...,L_{n}\}$\ there exist $L_{i_{j}}\in
\{L_{i_{1}},...,L_{i_{k}}\}$\ and $L_{j}\in \{L_{1},...,L_{n}\}$\ such that $%
lk(L_{j},L_{i_{j}})\neq 0$\ and $lk(L_{j},L_{i_{m}})=0$\ for each $%
L_{i_{m}}\in \{L_{i_{1}},...,L_{i_{k}}\} \backslash \{L_{i_{j}}\}.$
\end{corollary}

\begin{proof}
The previous property results from the proof of Lemma 4 in Chapter 9,
Paragraph I, p. 275 of \cite{Rolfsen} and more specifically from the form of
the simple closed curves $d_{i},$ $i=1,...,3g-1$\ on a closed surface $X$\
of genus $g,$ see Figure 1. In fact, it is a well-known result that any
orientation-preserving homeomorphism of a closed surface $X$ can be
expressed as a composition of Dehn twists along certain simple closed curves
belonging to the specified family $\{d_{i}\}.$\ For the curves $\{d_{i}\}$\
a similar property to Property (v) is valid. That is, for each family of
curves $\{d_{i_{1}},...,d_{i_{k}}\} \subseteq \{d_{1},...,d_{3g-1}\} $\
there exist $d_{i_{j}}\in $\ $\{d_{i_{1}},...,d_{i_{k}}\}$\ and $d_{j}\in
\{d_{1},...,d_{3g-1}\}$\ such that the intersection number $%
\#(d_{j},d_{i_{j}})=1$\ and $\#(d_{j},d_{i_{m}})=0$\ for each $d_{i_{m}}\in
\{d_{i_{1}},...,d_{i_{k}}\} \backslash \{d_{i_{j}}\}.$

Now, it also well known, that any closed 3-manifold $M$ can be obtained by
gluing two handlebodies $W,$ $W^{\prime }$ via an homeomorphism $h:\partial
W\rightarrow \partial W^{\prime }.$ Setting $X=\partial W,$ it results that $%
M$ is obtained by performing surgery on $S^{3}$ along tori which are the
boundaries of solid tori $V_{i}$ excavated from $W$ just under some of the $%
3g-1$ curves $d_{i}$ of Figure 1. We remark here that if we consider solid
tori which correspond to two parallel curves to $d_{i},$ say $d_{i}^{\prime
},$ $d_{i}^{\prime \prime },$ and if $T_{i}^{\prime },$ $T_{i}^{\prime
\prime }$ are the boundaries of the solid tori corresponding to $%
d_{i}^{\prime },$ $d_{i}^{\prime \prime },$ then there will be an annulus $A$
in $M$ connecting these tori i.e. one boundary component of $A$ will belong
to $T_{i}^{\prime }$ and the other will belong to $T_{i}^{\prime \prime }.$
This implies that $M$ is not acylindrical, as we have shown in Theorem \ref%
{link complement}. Therefore each curve $d_{i}$ generates exactly one solid
torus $V_{i}$ with core curve a component $L_{i}$ of the link $\mathcal{L}.$
So, the previous property concerning the curves $d_{i}$\ leads to Property
(v) which pertains to the components $L_{i}$ of $\mathcal{L}.$
\end{proof}

\begin{remark}
The fact that each component $L_{i}$ of $\mathcal{L}$ is unknotted is needed
to conclude that $\mathcal{L}$ consists at least of two components and that $%
N_{0}$ is acylindrical.
\end{remark}

%%%%%%%%%%%%%%%%%%%%%%%%%%%%%%%%%%%%%%%%%%
\begin{figure}[tbp]
%\hspace*{9mm}
\includegraphics[bb = 64 233 524 719]{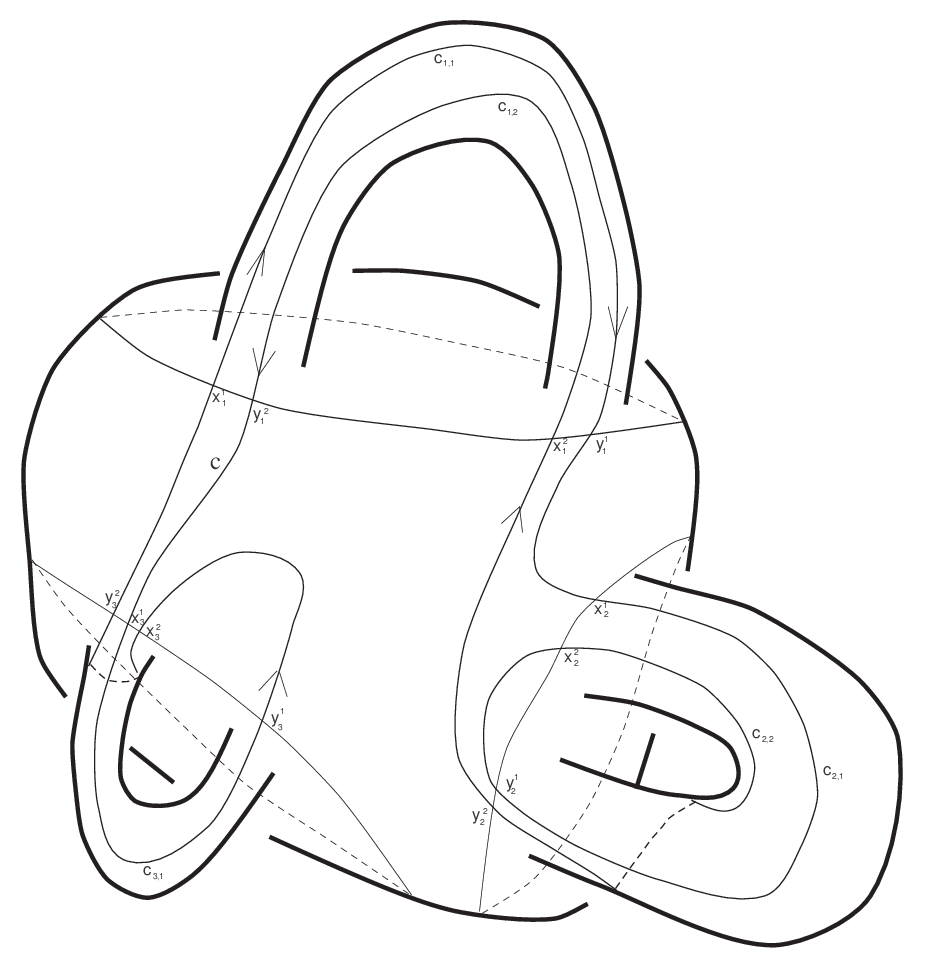} 
%\includegraphics[scale=0.94]{HAKEN_2.pdf} %\\vspace*{-8mm}\newline
%\label{disk4}
\caption{The curve $c$ which bounds a disc in $N.$}
\end{figure}

\begin{figure}[tbp]
%\hspace*{9mm}
\includegraphics[bb = 89 161 512 632]{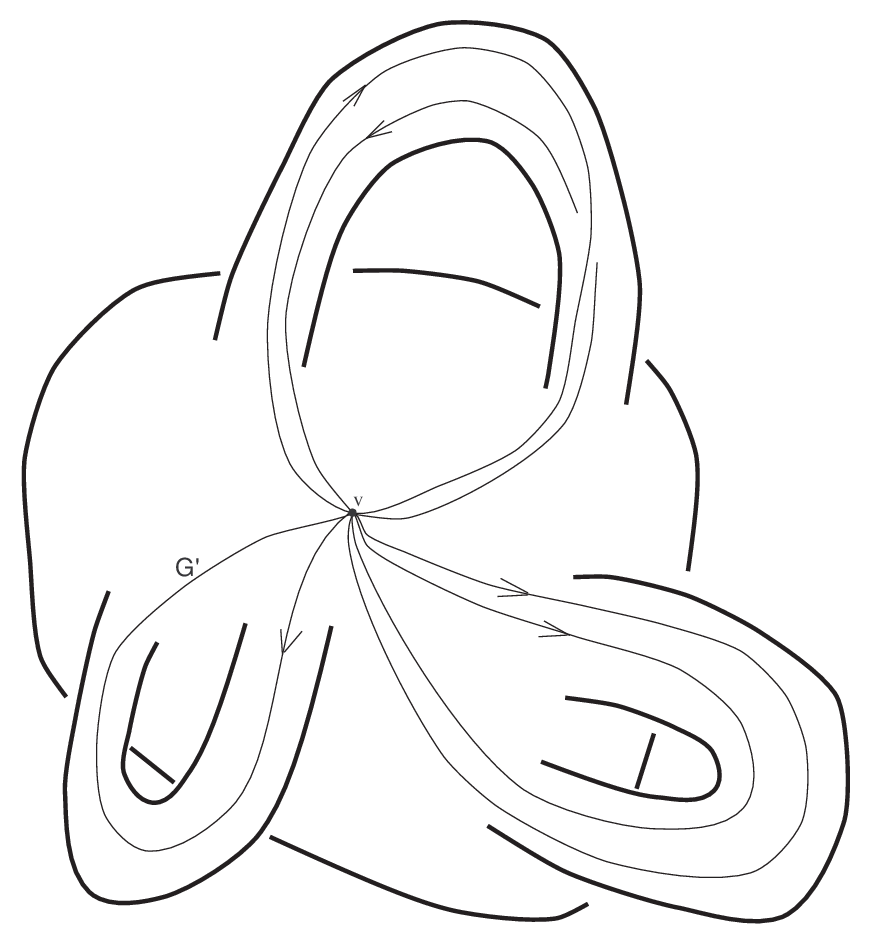} 
%\includegraphics[scale=0.94]{HAKEN_3.pdf} %\\vspace*{-8mm}\newline
%\label{disk4}
\caption{The graph $G^{\prime }$ that extends the spine $G$ of handlebody $%
H. $}
\end{figure}

%\begin{center}
\hspace*{-6cm} 
\begin{figure}[tbp]
%\hspace*{9mm}
\includegraphics[bb = 0 317 347 720]{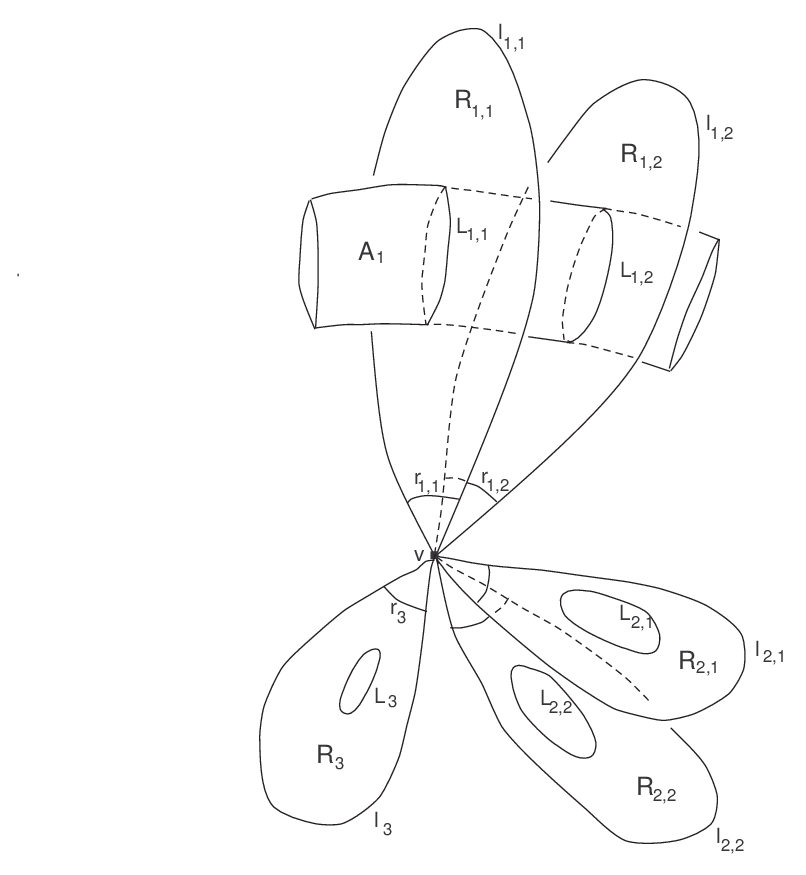} 
%\includegraphics[scale=0.94]{HAKEN_4.pdf} %\\vspace*{-8mm}\newline
%\label{disk4}
\caption{Step 1 for constructing the surface $Q.$}
\end{figure}
%\end{center}

\begin{figure}[tbp]
%\hspace*{9mm}
\includegraphics[scale=0.9, bb = 20 259 512 686]{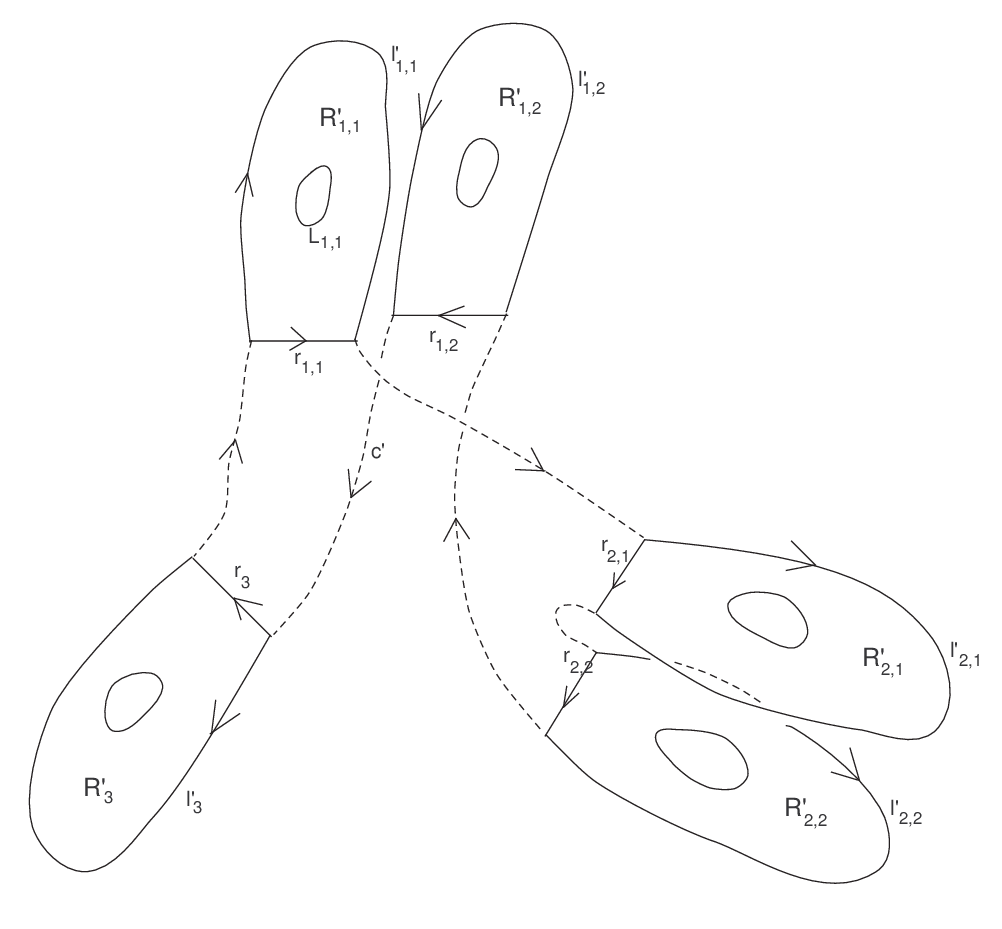} 
%\includegraphics[scale=0.94]{HAKEN_5.pdf} %\\vspace*{-8mm}\newline
%\label{disk4}
\caption{Step 2 for constructing the surface $Q.$}
\end{figure}

\begin{figure}[tbp]
%\hspace*{9mm}
\includegraphics[bb = 26 277 476 754]{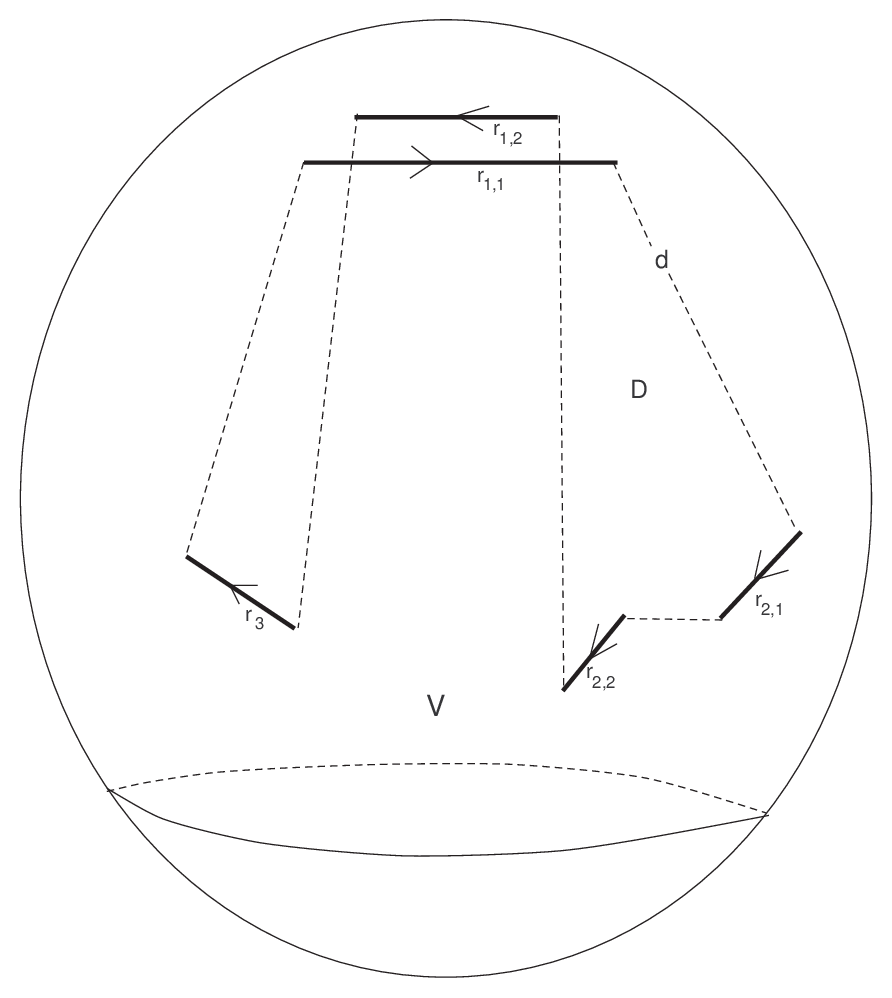} 
%\includegraphics[scale=0.94]{HAKEN_6.pdf} %\\vspace*{-8mm}\newline
%\label{disk4}
\caption{The curve $d$ in $V.$}
\end{figure}

\begin{figure}[tbp]
%\hspace*{9mm}
\includegraphics[bb = 67 352 428 627]{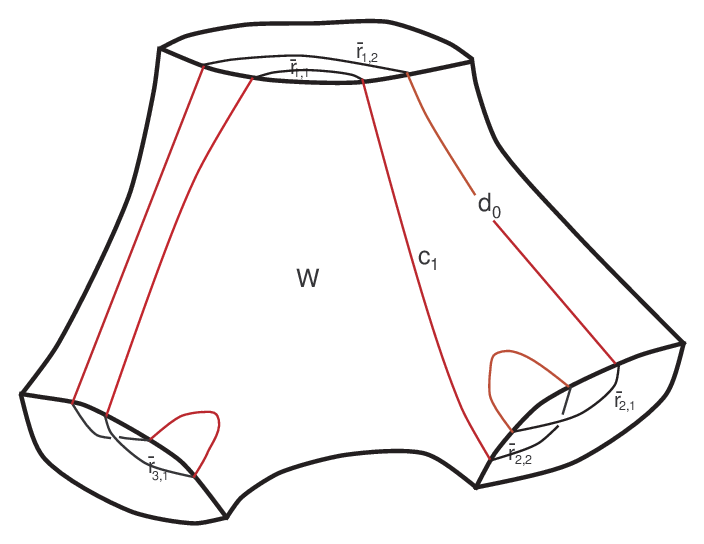} 
%\includegraphics[scale=0.94]{HAKEN_7.pdf} %\\vspace*{-8mm}\newline
%\label{disk4}
\caption{The curve $d_{0}$ in $W.$}
\end{figure}

\begin{figure}[tbp]
%\hspace*{9mm}
\includegraphics[bb = 26 298 448 629]{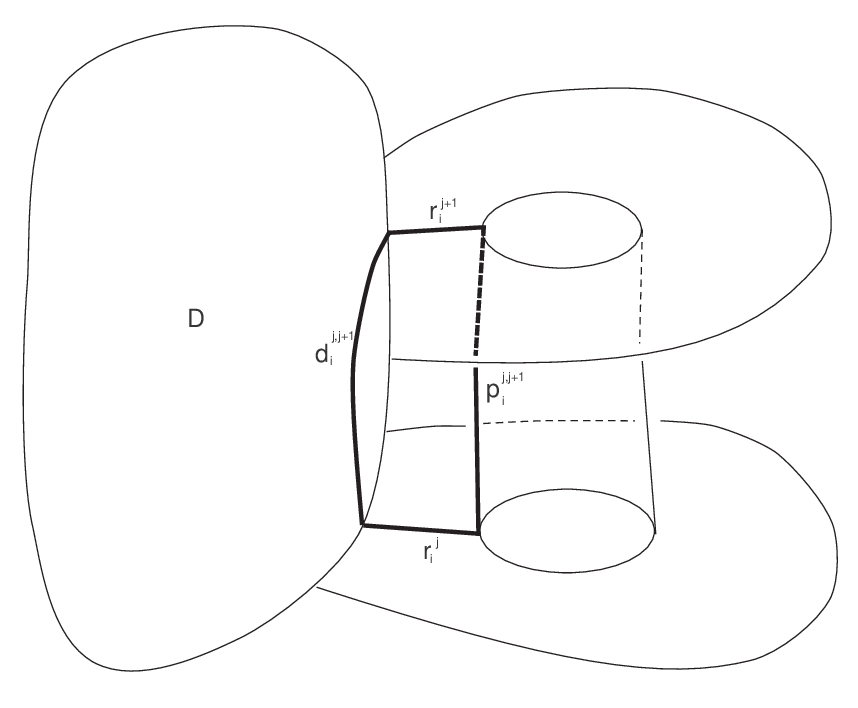} 
%\includegraphics[scale=0.64]{HAKEN_8.pdf} %\\vspace*{-8mm}\newline
%\label{disk4}
\caption{The construction of a meridian in $H_{0}.$}
\end{figure}
%%%%%%%%%%%%%%%%%%%%%%%%%%%%%%%%%%%%%%%%%%%%%%%%%%%%%%5
Now, we join the components $L_{i}$ and $L_{i+1}$ by a simple arc $a_{i}$
for each $i=1,...,n-1.$ Thus, we construct a graph $G_{0}$ in $S^{3}$ and we
consider a tubular neighborhood $H$ of $G_{0}.$ Obviously, $H$ is a
handlebody of genus $n$ and $\cup _{i}V_{i}\subset H$ therefore $%
S^{3}\backslash H\subset S^{3}\backslash (\cup _{i}V_{i}).$ Taking the
closure of $S^{3}\backslash (\cup _{i}V_{i})$ we get a compact manifold $%
N_{0}$ with $n$ toral boundary components $S_{i}.$ On each $S_{i}$ a simple
closed curve $s_{i}$ is fixed so that, if the solid torus $V_{i}$ is glued
to $S_{i}$ so that the meridian $m_{i}$ of $\partial V_{i}$ to match with $%
s_{i}$ we get the manifold $M.$ Now, let $N$ be the closure of $%
S^{3}\backslash H.$ In the boundary $\partial N$ the slopes $s_{i}$ continue
to live. It is clear that if we glue each meridian $m_{i}$ of $H$ to $s_{i}$
we obtain again the manifold $M.$ The manifold $N$ will be called the
complement of the handlebody $H$ corresponding to the link $\mathcal{L}$ and
will be referred to as the result of filling $N$ by that same handlebody $H.$

As as a consequence of Theorem \ref{link complement} we obtain the following
result, which is interesting in its own right.

\begin{theorem}
\label{attaching handlebody} Let $M$ be a closed orientable, irreducible,
atoroidal 3-manifold with infinite fundamental group. Then $M$ can be
recovered by gluing a handlebody $H$ of genus $g\geq 2$ in a manifold $N$
with boundary, via a homeomorphism $f:\partial N\rightarrow \partial H,$
such that $N$ is irreducible, atoroidal and boundary irreducible i.e. the
boundary $\partial N$ of $N$ is incompressible.
\end{theorem}

\begin{proof}
We use the terminology of Theorem \ref{link complement} and we recall that $%
V_{i}$ is a solid torus around $L_{i}$ for each $i$ and that $M$ is obtained
from $S^{3}$ by performing surgery along $\mathcal{L}=L_{1}\cup \cdots \cup
L_{n}.$ We will show that if join the components $L_{i},$ $L_{i+1}$ with
arbitrary arcs $a_{i},$ $i=1,..,n-1,$ and if we consider a tubular
neighborhood of the graph $G=L_{1}\cup a_{1}\cup L_{2}\cup a_{2}\cup \cdots
\cup a_{n-1}\cup L_{n}$ we obtain a handlebody $H$ such that $N=M\backslash
(Int(H))$ has the required properties.

More precisely, we have:

(1) $N$ \textit{is irreducible}.

Let $S$ be a 2-sphere in $N.$ Then $S$ is contained within the interior of $%
N_{0}.$ To demonstrate this, assume for contradiction that $S$ does not
separate $N_{0}.$ In this case, $S$ would not separate $M,$ implying that $%
M=M^{\prime }\#(S^{2}\times S^{1})$ for some closed 3-manifold $M^{\prime }.$
Now, let $C^{+}$ and $C^{-}$ be the two components of $N_{0}\backslash S.$
All solid tori $V_{i}$ must be contained within one of these components, say 
$C^{+}.$ If any $V_{i}$ were contained in $C^{-},$ then $S$ would not be
contained in $N$ as it would intersect the union of arcs $\cup
_{i=1}^{n-1}a_{i},$ contradicting the assumption that $S$ is contained in $%
N. $ Therefore, we conclude that $C^{-}$ is homeomorphic to a 3-ball as $S$
is embedded in $S^{3}$ and thus $C^{-}$ coincides with one component of $%
S^{3}\backslash S.$ Clearly, $C^{-}\subset N.$

(2) $N$ \textit{is atoroidal}.

Let $T$ be a torus in $N.$ First, we may prove that $T$ separates $N.$
Indeed, since $T\cap (\cup _{i=1}^{n-1}a_{i})\neq \varnothing $ we have that 
$T$ is contained in $N_{0}$ and separates it, as it is shown in the proof of
Theorem \ref{link complement}. Therefore, $T$ separates $N.$ Let $R^{+}$ and 
$R^{-}$ be the two components of $N\backslash T.$ Without loss of
generality, we may assume that all $V_{i}$ are contained in $R^{+}.$

From the proof of Theorem \ref{link complement}, we conclude that the torus $%
T$ bounds a solid torus $B\subset M$ on at least one side. If $B=R^{-}$ then 
$T$ is not incompressible in $N,$ that is, $N$ is atoroidal. Assuming $%
B=R^{+},$ we obtain a contradiction in a manner analogous to the proof of
Theorem \ref{link complement}.

(3) $\partial N$ is \textit{incompressible}.

We will illustrate this part of the proof with a concrete example, see
Figures 2 - 7.

First we need some terminology: A spine in a handlebody $H$ of genus $g\geq
1 $ is a graph $G$ embedded in $H$ so that $H\backslash G$ is homeomorphic
to $\partial H\times (0,1].$

In what follows we consider a spine $G\subset H$ such that $G$ consists of a
single vertex $v$ and $n$ loops $l_{1},...,l_{n},$ where each $l_{i}$ is
isotopic to the component $L_{i}$ of $\mathcal{L}$ within $H.$ Note that the
loops $l_{i},$ which comprise $G$ and the components $L_{i}$ of $\mathcal{L}%
, $ are considered unoriented. When we wish to emphasize that a loop (or
curve) $l$ is equipped with an orientation, we denote it by$\overrightarrow{l%
}$ (with an arrow on top).

Now suppose $c$ is a simple closed curve in $\partial N,$ as illustrated in
Figure 2, that bounds a disc $C$ embedded in $N.$ We will derive a
contradiction under this assumption. The curve $c$ on $\partial N$ can be
homotoped into the spine $G$ of the handlebody $H.$ Once in $G,$ this curve
decomposes into a concatenation of loops $l_{i_{1}},...,l_{i_{k}},$ where
each loop $l_{i_{j}}$ corresponds naturally to a component $L_{i_{j}}$ of
the link $\mathcal{L},$ for $j=1,...,k.$ Using the graph $G,$ we construct a
new graph $G^{\prime },$ which retains the same base vertex $v$ and depends
on the curve $c$ as follows: for each loop $l_{i_{m}}$ in $G,$ if the image
of $c$ winds around $l_{i_{m}}$ exactly $k_{i_{m}}$ times (ignoring
orientation), then we augment $G$ by adding $k_{i_{m}}-1$ additional loops
at the same handle containing $l_{i_{m}}.$ These added loops are called
parallel to $l_{i_{m}},$ and together they represent the equivalence class $%
[l_{i_{m}}];$ this process is performed for each $i_{m},$ $m=1,...,k,$ (see
Figure 3).

The curve $c,$ with a given orientation, will be denoted by $\overrightarrow{%
c}.$ For example the curve $\overrightarrow{c}$ of Figure 2 has the form 
\begin{equation*}
\overrightarrow{c}=\overrightarrow{l}_{1,1}\cdot \overrightarrow{l}%
_{2,1}\cdot \overrightarrow{l}_{2,2}\cdot \overrightarrow{l}_{1,2}\cdot 
\overrightarrow{l}_{3}.
\end{equation*}

Let $\mathcal{L}_{c}$ denote the sublink of $\mathcal{L}$ comprised of
components $L_{i_{1}},...,L_{i_{k}}.$ Define a new link $\mathcal{L}%
_{c}^{\prime }$ from $\mathcal{L}_{c}$ by expanding each component $%
L_{i_{m}} $ corresponding to a class $[l_{i_{m}}],$ which appears as $%
k_{i_{m}}$ parallel strands, into $k_{i_{m}}$ separate components in $%
\mathcal{L}_{c}^{\prime }.$ In other words, $\mathcal{L}_{c}^{\prime }$ lies
naturally above $G^{\prime },$ and compared to $\mathcal{L}_{c},$ $\mathcal{L%
}_{c}^{\prime }\backslash \mathcal{L}_{c}$ includes an extra $(k_{i_{m}}-1)$
copies of each component $L_{i_{m}}$ originally in $\mathcal{L}_{c}.$ All
copies corresponding to the same $L_{i_{m}}$ for a given $i_{m}$ are termed
parallel. By definition, there exists an annulus $A_{i_{m}}$ in $%
S^{3}\backslash (\cup _{i}V_{i})$ that contains all the copies of $L_{i_{m}} 
$ (see Figure 4). The knot class consisting of all knots parallel to $%
L_{i_{m}}$ is denoted by $[L_{i_{m}}].$

The existence of the disc $C$ in $N$ enables the construction of a genus-0
surface $Q$ in $S^{3}$ whose boundary components are the components of $%
\mathcal{L}_{c}^{\prime }.$ To construct $Q,$ we proceed with the following
steps:

\textbf{Step 1}: For each loop $l_{i}$ of $G$ that lies within a handle $%
H_{i}$ of $H$ we consider all parallel loops $%
l_{i,1},l_{i,2},...,l_{i,k_{i}} $ of $G^{\prime }$ in that same handle $%
H_{i};$ here, the index $(i,1)$ corresponds to the original index $i,$ i.e. $%
(i,1)\equiv i.$ We also consider the parallel knots $L_{i,1},...,L_{i,k_{i}}$
in $H_{i}$ and for each $j\in \{1,...,k_{i}\}$ we consider an annulus $%
R_{i,j}$ whose boundary satisfies $\partial R_{i,j}=l_{i,j}\cup L_{i,j}.$
Moreover, the intersection of all these annuli is precisely the single point 
$v;$ that is,%
\begin{equation*}
\cap _{j=1}^{k_{i}}R_{i,j}=\{v\}.
\end{equation*}%
Finally, for future reference, within each annulus $R_{i,j}$ we consider an
arc $r_{i,j}$ located in a neighborhood $V$ of the vertex $v$ (see Figure 4).

For example, consider the curve $c$ within the genus-3 handlebody shown in
Figure 2. We regard the knots $L_{1},L_{2}$ and $L_{3}$ as core loops in the
respective handles $H_{1},H_{2}$ and $H_{3}$ of $H.$ In particular, for each
of $L_{1}$ and $L_{2},$ we introduce two parallel knots, denoted by $%
L_{1,1}, $ $L_{1,2}$ and $L_{2,1},$ $L_{2,2}$ respectively. We then consider
the following annuli:

\begin{itemize}
\item $R_{1,1}$ with boundary $\partial R_{1,1}=l_{1,1}\cup L_{1,1},$

\item $R_{1,2}$ with boundary $\partial R_{1,2}=l_{1,2}\cup L_{1,2},$

\item $R_{2,2}$ with boundary $\partial R_{2,2}=l_{2,2}\cup L_{2,2},$

\item $R_{3}$ with boundary $\partial R_{3}=l_{3}\cup L_{3}.$
\end{itemize}

\noindent See Figure 4 for a visual depiction of these constructions.

\textbf{Step 2}: For each $i$ and $j,$ we focus on a sub-annulus $%
R_{i,j}^{\prime }$ of $R_{i,j}$ chosen so that the arc $r_{i,j}$ lies within
its boundary. We assign orientations to the arcs $r_{i,j}$ to match the
orientations induced by the loops $\overrightarrow{l}_{i,j}$ (see Figures 3
and 5). Consequently, the boundary of $R_{i,j}^{\prime }$ can be written as%
\begin{equation*}
\partial R_{i,j}^{\prime }=r_{i,j}\cup l_{i.j}^{\prime }\cup L_{i,j}
\end{equation*}%
where $l_{i.j}^{\prime }$ is a subsegment of $l_{i,j}$ (see Figure 5).

\textbf{Step 3}: We draw a neighborhood $V$ of the vertex $v$ within the
interior of $H;$ topologically, $V$ is homeomorphic to a closed 3-ball. The
interiors of all arcs $r_{i,j}$ lie within the interior of the 3-ball
neighborhood $V,$ while their endpoints lie on its boundary; formally,%
\begin{equation*}
Int(r_{i,j})\subset Int(V)~\text{and }r_{i,j}\cap \partial V=\partial
r_{i,j}.
\end{equation*}%
See Figure 6, where the arcs $r_{i,j}$ are shown in bold. In a collar
neighborhood of $\partial V,$ we insert arcs that connect all the endpoints
of the arcs $r_{i,j}.$ These new arcs together with the $r_{i,j}$ (within
the collar) form a single simple closed curve $d.$ This curve $d$ then
bounds a 2-dimensional disk $D$ embedded in $V.$ Crucially, there exists an
orientation on $d$ such that all the arcs $\overrightarrow{r}_{i,j}$ (that
is, the arcs $r_{i,j}$ with their given orientations), are coherently
aligned with this orientation on $d.$ Moreover, the sequence in which the $%
r_{i,j}$ are connected follows the path of the curve $\overrightarrow{c}.$
Consequently, $d$ is obtained by replacing each arc $l_{i,j}^{\prime },$
(from Figure 5) with the corresponding bold arc $r_{i,j}.$ For example, in
Figure 6, the curve $d$ is formed by linking the bold arcs $%
r_{1,1},r_{2,1},r_{2,2},r_{1,2},r_{3}$ using dotted connector arcs.

\textbf{Step 4}: We consider $n$ separating discs $D_{1},...,D_{n}$ in $H$
that separate the handles $H_{1},...,H_{n}$ of $H.$ Following the curve $%
\overrightarrow{c}$ i.e. the oriented curve $c,$ we consider the points of
intersection $c\cap (\cup _{i=1}^{n}d_{i}),$ see Figure 2. Let $x_{i}^{j}$
and $y_{i}^{j}$ denote the entry and exit points of the curve $c$ in the
handle $H_{i}$ along $d_{i}.$ These points correspond to the intersections
of $c$ with $\cup _{i=1}^{n}d_{i}.$ Even if the number of intersection
points $c\cap (\cup _{i=1}^{n}d_{i})$ can be reduced, we assume that each
time $c$ enters the handle $H_{i},$ it intersects $d_{i}$ at exactly two
points: one upon entering and one upon leaving. A typical example of this
situation is the arc $c_{2,2}$ in Figure 2. We denote by $c_{i,j}$\ the
subarc of $c$\ in the handle $H_{i}$ with boundary\textbf{\ }$\partial
c_{i,j}=\{x_{i}^{j},y_{i}^{j}\}.$

Let $W$ be the subset of $H$ such that $\cup _{i=1}^{n}D_{i}\subset \partial
W,$ where each $D_{i}$ is a separating disk of $H.$ In this context, $W$ can
be viewed as a spotted 3-ball. Consider arcs $\overline{r}_{i,j}$ such that $%
Int(\overline{r}_{i,j})\subset Int(W)$ and $\partial \overline{r}%
_{i,j}=\{x_{i}^{j},y_{i}^{j}\}.$ Additionally, let $c_{i},$ for $i=1,...,a,$
be the subarcs of $c\cap \partial W.$ All arcs $c_{i}$ are drawn in red in
Figure 7, and among them the arc $c_{1}$ is highlighted. Then the union of
all $\overline{r}_{i,j}$ and all $c_{i}$ forms a simple closed curve $d_{0}.$
This curve bounds an embedded disk $D_{0}\subset W$ (see Figures 2 and 7).
Now, it is clear that we can construct an ambient isotopy $f_{t},$ $t\in
\lbrack 0,1],$ that carries $d_{0}$ to $d,$ mapping each arc $\overline{r}%
_{i,j}$ to $r_{i,j}$ while avoiding the disc $D.$ Furthermore, this isotopy
can be extended so that each boundary subarc $c_{i,j}$ is carried to the arc 
$l_{i,j}^{\prime }.$

\textbf{Step 5}: After constructing the disc $D$ (see Figure 6), we glue to
it the annuli $R_{i,j}^{\prime }$ illustrated in Figure 5 along the arcs $%
r_{i,j}.$ Thus, we produce a planar surface $Q_{1}$ bounded by a curve $%
c^{\prime }$ and the curves $L_{i,j}.$ Thanks to the isotopy $f_{t}$ from
Step 4, the curve $c^{\prime }$ is freely homotopic to $c$ within the
handlebody $H.$ Therefore, there exists an annulus $C^{\prime }\subset H$
with boundary $\partial C^{\prime }=c^{\prime }\cup c$ and we may assume
that $C^{\prime }\cap Q_{1}=c^{\prime }.$ Gluing $C^{\prime }$ to $Q_{1}$
along $c^{\prime }$ yields a new surface $Q_{2}.$ Next, we glue the
compressing disc $C\subset N$ with $Q_{2}$ along the curve $c.$ The
resulting surface is the desired surface $Q.$

Once the surface $Q$ has been constructed, we distinguish the following two
cases (i) and (ii):

(i) $\mathcal{L}_{c}=\mathcal{L}_{c}^{\prime }.$

In this case, Property (v) implies that there is a component, say $%
L_{i_{0}}, $ of $\mathcal{L}$ which has linking number $\pm 1$ with some
component, say $L_{i_{0}}^{\prime },$ of $\mathcal{L}_{c}$ and $L_{i_{0}}$
has linking number $0$ with all the other components of $\mathcal{L}_{c}.$
We denote by $\Delta _{i_{0}}$ and $\Delta _{i_{0}}^{\prime }$ embedded
2-discs in $S^{3}$ with $\partial \Delta _{i_{0}}=L_{i_{0}}$ and $\partial
\Delta _{i_{0}}^{\prime }=L_{i_{0}}^{\prime }.$ Since $%
lk(L_{i_{0}},L_{i_{0}}^{\prime })=\pm 1$ it follows that $\Delta
_{i_{0}}\cap \Delta _{i_{0}}^{\prime }$ is a single arc $a$ whose endpoints $%
\partial a=\{p,q\}$ satisfy $p\in L_{i_{0}}$ and $q\in L_{i_{0}}^{\prime }.$
This fact, together with Property (v), implies that $\Delta _{i_{0}}\cap Q$
contains an arc starting at $q$ which cannot end on any other component of
component of $\mathcal{L}_{c}$ except $L_{i_{0}}^{\prime }.$ Therefore, $%
\Delta _{i_{0}}\cap L_{i_{0}}^{\prime }$ must consist of more than one
point, which contradicts the assumption that $\partial \Delta
_{i_{0}}^{\prime }=L_{i_{0}}^{\prime }$ meets $\Delta _{i_{0}}$ exactly in
the single point $q.$

(ii)\textbf{\ }$\mathcal{L}_{c}\neq \mathcal{L}_{c}^{\prime }.$

In this case, we distinguish below two subcases, (ii$_{a})$ and (ii$_{b}),$
which will be considered separately.

$(ii_{a})$\textbf{\ }\textit{For each }$L_{i_{m}}\in \mathcal{L}_{c}^{\prime
},$\textit{\ the class }$[L_{i_{m}}]$\textit{\ contains an even number of
copies parallel to }$L_{i_{m}}.$

Let $L_{i_{m},j},$ $j=1,...,2k_{i_{m}}$ be the elements of $[L_{i_{m}}].$ We
may assume that there exists an annulus $A_{i_{m}}$ such that%
\begin{equation*}
A_{i_{m}}\cap Q=\{L_{i_{m},j}:j=1,...,2k_{i_{m}}\}
\end{equation*}%
and if we identify $A_{i_{m}}$ with $S^{1}\times \lbrack 0,1],$ then the
copies $L_{i_{m},j}$ correspond to the circles $S^{1}\times \{t_{j}\},$ $%
j=1,...,2k_{j_{m}},$ where $t_{1}<\cdots <t_{2k_{i_{m}}};$ such an annulus,
for example the annulus $A_{1},$ is illustrated in Figure 4. Recall that,
for fixed $i_{m},$ the annuli $R_{i_{m},j}^{\prime }$ and $%
R_{i_{m},j+1}^{\prime }$ can be considered parallel (see Figure 5) in the
sense: there is a neighborhood $U$ in $H$ homeomorphic to $S^{1}\times
\lbrack 0,1]\times \lbrack 0,1]$ such that

\begin{itemize}
\item $S^{1}\times \lbrack 0,1]\times \{0\}$ corresponds to $%
R_{i_{m},j}^{\prime },$

\item $S^{1}\times \lbrack 0,1]\times \{1\})$ corresponds to $%
R_{i_{m},j+1}^{\prime },$

\item $Int(U)\cap R_{i_{m},j}^{\prime }=\emptyset $ for any $%
j=1,...,2k_{i_{m}}.$
\end{itemize}

Now we form an orientable closed surface $P$ in $H$ as follows: for each $%
i_{m}\in \{i_{1},...,i_{k}\}$ and each $j\in \{1,3,...,2k_{i_{m}}-1\},$ let $%
A_{i_{m}}^{j,j+1}\subset A_{i_{m}}$ be the annulus with $\partial
A_{i_{m}}^{j,j+1}=L_{i_{m},j}\cup L_{i_{m},j+1}.$ By construction, these
annuli are disjoint. We glue $A_{i_{m}}^{j,j+1}$ to $Q$ along the common
boundary components. The resulting surface is the surface $P.$ By
construction (and up to isotopy) we have that $L\cap P=\emptyset $ for each $%
L\in $ $\mathcal{L}\backslash \mathcal{L}_{c}.$

We now have the following claim.

\textit{Claim 1}: The surface $P$ bounds a handlebody $H_{0}.$

\textit{Proof of Claim 1}. To prove that $P$ bounds a handlebody, it
suffices to exhibit a collection of disjoint, non-separating compressing
discs of $P$ along which $P$ is cut to a 2-sphere. We construct such
compressing discs:

\begin{itemize}
\item For each $i=i_{1},...,i_{k}$ and for each $j=1,3,...,2k_{i_{m}}-1,$ we
define a disc $D_{i}^{j,j+1}$ whose boundary $\partial D_{i}^{j,j+1}$
consists of an arc $p_{i}^{j,j+1}\subset A_{i}^{j,j+1},$ two arcs $%
r_{i}^{j}\subset R_{i,j}^{\prime }$ and $r_{i}^{j+1}\subset
R_{i,j+1}^{\prime }$ and an arc $d_{i}^{j,j+1}\subset D,$ see Figure 8.
\end{itemize}

The intersection of two such compressing discs is empty. Specifically, the
boundaries of these discs, denoted as $D_{1}$ and $D_{2}$ for simplicity,
can intersect only along the arcs $\partial D_{1}\cap D$ and $\partial
D_{2}\cap D.$ However, these arcs can intersect at most once, up to isotopy.
This implies that $\partial D_{1},\partial D_{2}$ intersect at most once.
Such an intersection is impossible, as two compressing discs of $P$ must
intersect an even number of times. Therefore $H_{0}$ must be a handlebody.
Finally, by construction, and up to isotopy, there are no components of $%
\mathcal{L}$\ in the interior of $H_{0}.$ This completes the proof of Claim
1.

Next, we have

\textit{Claim 2}: Let $H_{0}^{c}$ be the closure of the complement of $H_{0}$
in $S^{3}.$ Then for each $i=i_{1},...,i_{k}$ and each $j=1,...,k_{i_{m}},$
we may find disjoint discs $\Delta _{i,j}$ with $\partial \Delta
_{i,j}=L_{i,j}$ and $\Delta _{i,j}\subset H_{0}^{c}.$ Therefore $H_{0}^{c}$
is a handlebody.

\textit{Proof of Claim 2}. From the construction of $H_{0}$ we may find
discs $\Delta _{i,j}$ with $\partial \Delta _{i,j}=L_{i,j}$ and $\Delta
_{i,j}\cap H_{0}=L_{i,j}$ for each $i=i_{1},...,i_{k}$ and $%
j=1,...,k_{i_{m}}.$ Since the components $L_{i,j}$ are disjoint for all $%
i,j, $ the discs $\Delta _{i,j}$ must also be disjoint. In fact, suppose,
for contradiction, that two discs $\Delta _{i_{1},j_{1}}$ and $\Delta
_{i_{2},j_{2}}$ intersect. Let $x=L_{i_{2},j_{2}}\cap Int(\Delta
_{i_{1},j_{1}}).$ Then a neighborhood $U$ of $x$ in $\Delta _{i_{1},j_{1}}$
must satisfy $U\cap H_{0}=\emptyset .$ This is a contradiction because $x\in
\partial H_{0},$ implying that $\Delta _{i_{1},j_{1}}$ and $\Delta
_{i_{2},j_{2}}$ cannot intersect. Therefore, the discs $\Delta _{i,j}$ are
disjoint, and Claim 2 is proven.

Recall that $\mathcal{L}_{c}=\{L_{i_{1}},...,L_{i_{k}}\}.$ Then it is
important to observe that in that in situation $(ii_{a}),$ $\mathcal{L}%
\backslash \mathcal{L}_{c}\neq \emptyset .$ Indeed, if $\mathcal{L}%
\backslash \mathcal{L}_{c}=\emptyset ,$ then the construction of the
handlebody $H_{0}^{c}$ would imply that all the elements of $\mathcal{L}$
have linking number $0$ with one another, which is impossible.

\textit{Claim 3}: All the elements of $\mathcal{L}$ are contained in the
complement $H_{0}^{c}$ of $H_{0}.$

\textit{Proof of Claim 3}. There is an element $L\in \mathcal{L}\backslash 
\mathcal{L}_{c}$ which has linking number $\pm 1$ with some element of $%
\mathcal{L}_{c}.$ This implies that, up to isotopy, $L\subset
Int(H_{0}^{c}). $ Define $\mathcal{L}_{c}^{1}=\mathcal{L}_{b}\cup \{L\}.$
Then there is an element $L^{\prime }\in \mathcal{L}\backslash \mathcal{L}%
_{b}^{1}$ which has linking number $\pm 1$ with some element of $\mathcal{L}%
_{c}^{1}.$ Again, up to isotopy it follows that $L^{\prime }\subset
Int(H_{0}^{c}).$ Continuing in this way, one builds a sequence of elements
in $\mathcal{L}\backslash \mathcal{L}_{c},$ each of which is isotoped into
the interior of $H_{0}^{c}.$ This process proves Claim 3.

To reach a contradiction, we may also assume, up to isotopy, that all the
elements of $\mathcal{L}$ are contained in the interior of $H_{0}^{c}.$
Therefore we may find a surface $S$ isotopic to $\partial H$ lying in the
interior of $H_{0}^{c}.$ However, by construction the surfaces $P$ and
surface $Q$ intersect transversely along the simple closed curve $c.$ This
contradiction (namely, the existence of the surface $P$ with those
properties) shows that case $(ii_{a})$ cannot occur.

$(ii_{b})$\textbf{\ }\textit{There is at least one }$L_{i_{m}}\in \mathcal{L}%
_{c}^{\prime },$\textit{\ such that }$[L_{i_{m}}]$\textit{\ contains an odd
number of copies parallel to }$L_{i_{m}}.$

We approach this case by integrating the previous ones. More concretely, as
in Case (ii$_{a})$ we can construct a compact surface $P$ with boundary such
that each boundary component of $P$ corresponds to an element of some class $%
[L_{i_{m}}],$ and no two boundary components belong to the same class. This
assertion follows from the fact that each boundary component of $P$ comes
from each class $[L_{i_{m}}]$ which contains an odd number of elements.

Let $L_{b,1},...,L_{b,i}$ be the boundary components of $P$ and let $%
\mathcal{L}_{b}=\{L_{b,1},...,L_{b,i}\} \subset \mathcal{L}.$ By applying
Property (v) for the set $\mathcal{L}_{b}$ we can find an element $%
L_{i_{0}}\in \mathcal{L}$ and an element $L_{b,i_{0}}\in \mathcal{L}_{b}$
such that $lk(L_{i_{0}},L_{b,i_{0}})=1$ and $lk(L_{i_{0}},L_{b,i})=0$ for
each $L_{b,i}\in \mathcal{L}_{b}\backslash \{L_{b,i_{0}}\}.$ Let $\Delta
_{i_{0}}$ be a disc with $\partial \Delta _{i_{0}}=L_{i_{0}}$ and consider
the intersection $\Delta _{i_{0}}\cap P.$ This intersection set must contain
an arc starting at a point on the component $L_{b,i_{0}}$ of $P$ and which
cannot extend to another component of $P.$ Therefore, we reach a
contradiction as in Case (i).
\end{proof}

\begin{corollary}
\label{WLG}The interior of manifold $N$ of Corollary \ref{attaching
handlebody} admits a complete hyperbolic structure.
\end{corollary}

\begin{proof}
The result follows from Theorem \ref{Hyperbolization for the interior}.
\end{proof}

\begin{remark}
An interesting question is whether one can choose the arcs $a_{i},$ which
join the components $L_{i}$ of $\mathcal{L},$ so that $N$ becomes
acylindrical. If this were true, by Theorem \ref{Hyperbolization for compact}%
, $N$ would admit a hyperbolic structure with geodesic boundary. However,
the main result of this article does not hinge on this question, and so we
do not pursue it further.
\end{remark}

\section{A special surface in the complement of the handlebody corresponding
to link $\mathcal{L}$}

We begin by establishing the terminology. Let $N$ be an orientable, compact,
connected 3-manifold with connected incompressible boundary which is a
surface of genus $g\geq 2.$ By selecting disjoint non-separating curves $%
s_{i},$ $i=1,...,g$ on $\partial N,$ which divide $\partial N$ into a
2-sphere with $2g$ holes, we may construct a closed manifold $\overline{N}$
as follows. If $D^{2}$ denotes the 2-ball, we attach a copy of $D^{2}\times
I $ to $N$ by identifying $\partial D^{2}\times I$ with a neighborhood of $%
s_{i}$ in $\partial N,$ for each $i.$ The resulting manifold $N^{\prime }$
has a 2-sphere boundary. We obtain $\overline{N}$ by attaching a 3-ball $B$
to $N^{\prime }$ by identifying $\partial B$ with $\partial N^{\prime }.$
The curves $s_{i}$ will be referred to as the \textit{slopes} of $\partial N$
and we will say that $\overline{N}$ is obtained by gluing a handlebody $H$
of genus $g$ to $N$ along $\partial N$ via a homeomorphism $f:\partial
N\rightarrow \partial H$ so that the slopes $s_{i}$ correspond to disjoint
non-separating meridians of $H.$ The manifold $\overline{N}$ depends solely
on the slopes $s_{i}$ and is termed the \textit{filling} of\textbf{\ }$N$ by
a handlebody.

We begin with a compact, orientable, atoroidal, irreducible
three-dimensional manifold $M$ with an infinite fundamental group. As
discussed in Section 4, $M$ is constructed by considering a special link $%
\mathcal{L}$ in $S^{3},$ removing solid tori $V_{i}$ along each component $%
L_{i}$ of $\mathcal{L}$ and then gluing back $V_{i}$ in a different manner.
By Corollary \ref{attaching handlebody}, $M$ can also be viewed as the
filling of a manifold $N$ by a handlebody, where $N$ is orientable, compact,
atoroidal and has connected incompressible boundary. Recall that $%
N_{0}=S^{3}\backslash (\cup _{i}Int(V_{i}))$ and let $S_{i}$ be the boundary
components of $N_{0}$ such that $s_{i}\subset S_{i}$ for each $i.$

Our strategy for proving the Haken Conjecture can be outlined as follows.
Assuming that $M$ is non-Haken, we construct, using ideas from \cite%
{Cooper-Long-Reid}, an essential surface $S$ in $N$ with a specific property
described in the proposition below. In the subsequent section, the surface $%
S $ will be lifted to an embedded surface $\widetilde{S}$ in a finite cover $%
\widetilde{N}$ of $N.$ Then by compressing $\widetilde{S}$ in the filling
manifold $\widetilde{M}$ of $\widetilde{N}$ and projecting the result back
to $M,$ we obtain our desired outcome.

The objective of this section is to prove the following fundamental
proposition:

\begin{proposition}
\label{nice incompressible}If $M$ is not Haken or if $M$ does not have a
finite covering which is a surface bundle over $S^{1}$ then we may find a
properly embedded incompressible surface $T$ in $N$ such that $T\cap 
\mathcal{S}=\emptyset ,$ where $\mathcal{S}=\{s_{i},$ $i=1,\cdots ,g\}.$
\end{proposition}

\begin{proof}
The proof will be based on several claims, each of which is established
below.

We fix a set of simple closed curves $\{c_{1},\cdots ,c_{g},b_{1},\cdots
,b_{g}\}$ on $\partial N$ such that: $c_{i}=s_{i}$ for each $i=1,...,g$ and $%
b_{i}$ is the meridian (up to isotopy) of $V_{i}.$ We have that $b_{i}\cap
c_{i}$ is a single point for each $i$ and $b_{i}\cap c_{j}=\emptyset $ for $%
j\neq i.$

The homology classes of $c_{i}$ and $b_{i}$ on $\partial N$ generate the
homology group $H_{1}(\partial N)=\mathbb{Z}^{2g}.$ Let us denote by $%
[c_{i}],$ $[b_{i}]$ the homology class of $c_{i}$ and $b_{i}$ in $N.$

\textit{Claim 1}: Each homology class $[c_{i}]$ and $[b_{i}]$ has infinite
order in $H_{1}(N).$

\textit{Proof of Claim 1}. Recall that $V_{i}$ is an open tubular
neighborhood of the component $L_{i}$ of the link $\mathcal{L}$ in $S^{3}.$
Therefore, we may assume that $V_{i}\subset \cup _{i}V_{i}\subset H$ and
hence $N=S^{3}\backslash H\subset S^{3}\backslash (\cup _{i}V_{i})\subset
S^{3}\backslash V_{i},$ $\forall i.$ Now if $[c_{i}]$ (resp. $[b_{i}])$ has
finite order in $H_{1}(N)$ then there is a surface (incompressible) $%
F_{i}\subset N$ such that $\partial F_{i}$ consists of $k$ copies of the
curve $c_{i},$ $k\geq 1$ (see for instance, Proposition 2.6 in \cite%
{Fomenko-Matveev}). Obviously this surface $F_{i}$ lives in $%
N_{0}=S^{3}\backslash (\cup _{i}V_{i})$ and more specifically in $%
S^{3}\backslash V_{i}.$ Without loss of generality, we may assume that $%
F_{i} $ is incompressible in $S^{3}\backslash V_{i}.$ In fact, if $F_{i}$ is
compressible, then by compressing it within $S^{3}\backslash V_{i}$ we
obtain an incompressible component, say $F_{i}^{\prime },$ in $%
S^{3}\backslash V_{i}.$

Now, from the properties of the link $\mathcal{L},$ $S^{3}\backslash V_{i}$
is a solid torus. Therefore, we obtain a contradiction to the existence of
such a surface $F_{i}^{\prime }$ in $S^{3}\backslash V_{i}$ since the only
properly embedded incompressible surfaces in a solid torus are either
meridians or annuli parallel to the boundary. Consequently, the only case we
need to exclude is that $c_{i}$ or $b_{i}$ cannot be the boundary of a
properly embedded disc in $S^{3}\backslash V_{i}$ that is non-parallel to $%
\partial V_{i}.$ This arises from the fact that $c_{i}$ (resp. $b_{i})$ is a
curve of $\partial V_{i}$ which, when viewed in $V_{i},$ goes once around
the meridian and once around the longitude of $V_{i}$ (resp. $b_{i}$ bounds
a meridian of $V_{i}).$ Thus, Claim 1 is proven.

Let $G=H_{1}(N)$ and let $K$ be the image of $H_{1}(\partial N)$ by the
natural homomorphism $\tau _{\#}:$ $H_{1}(\partial N)\rightarrow H_{1}(N).$
So $K$ is a subgroup of $G.$

\textit{Claim 2}: $G$ and $K$ are finitely generated abelian groups and the
number of generators of infinite order of $G$ is greater than or equal to
the number of generators of infinite order of $K.$

\textit{Proof of Claim 2}. Since the fundamental group $\pi _{1}(N)$ of $N$
is a finitely presented group, $H_{1}(N)$ is finitely generated. Let $%
c_{1},\cdots ,c_{n}$ be the generators of $G.$ The group $qG$ is defined to
be the finitely generated abelian subgroup of $G$ with generators $%
qc_{1},\cdots ,qc_{n}.$ Since $K\leq G$ we have that $qK\leq qG$ for each $%
q\in \mathbb{Z}.$ From the structural theorem of finitely generated abelian
groups there is $p\in \mathbb{Z}$ such that $pG$ consists of elements of
infinite order i.e. it is a free abelian group. Obviously, an element $g$ of 
$G$ has infinite order in $G$ if and only if $pg$ has infinite order in $pG.$
Therefore, since $pK\leq pG$ the group $pK$ is free and the rank of the free
abelian group $pG$ is greater than or equal to the rank of $pK.$ Thus, Claim
2 is proven.

\textit{Case I}: The map $\tau _{\#}:$ $H_{1}(\partial N)\rightarrow
H_{1}(N) $ is a monomorphism.

\textit{Case II}: The map $\tau _{\#}:$ $H_{1}(\partial N)\rightarrow
H_{1}(N)$ is not a monomorphism.

Let $\tau :$ $\pi _{1}(N,x_{0})\rightarrow H_{1}(N)$ be the natural
epimorphism which connects the first fundamental group with the first
homology group of $N.$ For each element $[a]\in H_{1}(N)$ let us denote by $%
\overline{a}$ each element of $\tau ^{-1}([a]).$ Finally we remark that the
natural homomorphism $\iota :\pi _{1}(\partial N,x_{0})\rightarrow \pi
_{1}(N,x_{0})$ is $1-1$ since $\partial N$ is incompressible.

\textit{Claim 3}: In Case I, we may find an epimorphism $\phi :$ $\pi
_{1}(N)\rightarrow \mathbb{Z}$ such that each $\overline{c_{i}}$ is send to $%
0$ in $\mathbb{Z}$ but at least one $\overline{b_{i}}$ is send to a
non-trivial element of $\mathbb{Z}.$

\textit{Proof of Claim 3}. We first construct an epimorphism $\psi
:H_{1}(N)\rightarrow $ $\mathbb{Z}$ which has similar properties with $\phi $
and compose $\psi $ with the natural epimorphism $\tau :$ $\pi
_{1}(N,x_{0})\rightarrow H_{1}(N).$ Let $r_{j},$ $j=1,\cdots ,k$ be the
generators of $G$ of infinite order and $r_{j}^{\prime },$ $j=1,\cdots ,m$
be the generators of $G$ of finite order. From Claim 2 we have that $k\geq
2g.$ Each $[c_{i}]$ and each $[b_{j}]$ can be expressed as a linear
combination with integer coefficients of $r_{i}$ and of $r_{j}^{\prime }.$
So, we assume that we have

$[c_{1}]=\alpha _{11}r_{1}+\cdots +\alpha _{1k}r_{k}+\alpha _{11}^{\prime
}r_{1}^{\prime }+\cdots +\alpha _{1m}^{\prime }r_{m}^{\prime }$

$...................................................$

$[c_{g}]=\alpha _{g1}r_{1}+\cdots +\alpha _{gk}r_{k}+\alpha _{g1}^{\prime
}r_{1}^{\prime }+\cdots +\alpha _{gm}^{\prime }r_{m}^{\prime }$

$[b_{1}]=\beta _{11}r_{1}+\cdots +\beta _{1k}r_{k}+\beta _{11}^{\prime
}r_{1}^{\prime }+\cdots +\beta _{1m}^{\prime }r_{m}^{\prime }$

$...................................................$

$[b_{g}]=\beta _{g1}r_{1}+\cdots +\beta _{gk}r_{k}+\beta _{g1}^{\prime
}r_{1}^{\prime }+\cdots +\beta _{gm}^{\prime }r_{m}^{\prime }$

Our objective is to construct an epimorphism $\psi :H_{1}(N)\rightarrow $ $%
\mathbb{Z}$ which sends each $[c_{i}]$ and each $[b_{i}]$ to $0$ in $\mathbb{%
Z},$ for all $i\neq g$ and sends $[b_{g}]$ to a non-trivial element of $%
\mathbb{Z}.$ To achieve this, we associate a variable $x_{i}$ to each $r_{i}$
and a variable $x_{i}^{\prime }$ to each $r_{i}^{\prime }.$ We then replace
the previous equations with the following linear, non-homogeneous system:

$0=\alpha _{11}x_{1}+\cdots +\alpha _{1k}x_{k}+\alpha _{11}^{\prime
}x_{1}^{\prime }+\cdots +\alpha _{1m}^{\prime }x_{m}^{\prime }$

$...................................................$

$0=\alpha _{g1}x_{1}+\cdots +\alpha _{gk}x_{k}+\alpha _{g1}^{\prime
}x_{1}^{\prime }+\cdots +\alpha _{gm}^{\prime }x_{m}^{\prime }$

$0=\beta _{11}x_{1}+\cdots +\beta _{1k}x_{k}+\beta _{11}^{\prime
}x_{1}^{\prime }+\cdots +\beta _{1m}^{\prime }x_{m}^{\prime }$

$...................................................$

$1=\beta _{g1}x_{1}+\cdots +\beta _{gk}x_{k}+\beta _{g1}^{\prime
}x_{1}^{\prime }+\cdots +\beta _{gm}^{\prime }x_{m}^{\prime }$

The above (non-homogeneous) system has a solution, say $X=(\rho _{1},\cdots
,\rho _{k},\rho _{1}^{\prime },\cdots ,\rho _{m}^{\prime }).$ In fact,
because $k\geq 2g,$ it suffices to prove that the rank of the linear
homogeneous system, corresponding to the considered non-homogeneous system
is equal to $2g.$ This follows from the fact that the map $\tau _{\#}:$ $%
H_{1}(\partial N)\rightarrow H_{1}(N)$ is a monomorphism, which implies that
the matrix formed by the coefficients of the non-homogeneous system i.e the
matrix 
\begin{equation*}
\left( 
\begin{array}{cccccc}
\alpha _{11} & \cdots & \alpha _{1k} & \alpha _{11}^{\prime } & \cdots & 
\alpha _{1m}^{\prime } \\ 
\cdots & \cdots & \cdots & \cdots & \cdots & \cdots \\ 
\alpha _{g1} & \cdots & \alpha _{gk} & \alpha _{g1}^{\prime } & \cdots & 
\alpha _{gm}^{\prime } \\ 
\beta _{11} & \cdots & \beta _{1k} & \beta _{11}^{\prime } & \cdots & \beta
_{1m}^{\prime } \\ 
\cdots & \cdots & \cdots & \cdots & \cdots & \cdots \\ 
\beta _{g1} & \cdots & \beta _{fk} & \beta _{g1}^{\prime } & \cdots & \beta
_{gm}^{\prime }%
\end{array}%
\right)
\end{equation*}%
has rank $2g.$

By multiplying the solution $X=(\rho _{1},\cdots ,\rho _{k},\rho
_{1}^{\prime },\cdots ,\rho _{m}^{\prime })$ by an integer $l,$ we may
assume that the coordinates of $lX=(l\rho _{1},\cdots ,l\rho _{k},l\rho
_{1}^{\prime },\cdots ,l\rho _{m}^{\prime })$ are integers. Additionally, by
multiplying all the previous equations by $l,$ we obtain a new system which
has as a solution $lX=(l\rho _{1},\cdots ,l\rho _{k},l\rho _{1}^{\prime
},\cdots ,l\rho _{m}^{\prime }).$

Now we are ready to define the epimorphism $\psi :H_{1}(N)\rightarrow $ $%
\mathbb{Z}$ as follows: It suffices to define $\psi $ on the set of
generators of $H_{1}(N).$ Thus we set $\psi (r_{1})=l\rho _{1},\cdots ,\psi
(r_{k})=l\rho _{k},\psi (r_{1}^{\prime })=l\rho _{1}^{\prime },\cdots ,\psi
(r_{m}^{\prime })=l\rho _{m}^{\prime }.$ Then, the previous equations assert
that $[b_{g}]$ is send to $l\in \mathbb{Z}$ while all the other $[b_{i}]$
and $[c_{i}]$ are sent to $0\in \mathbb{Z}.$

Finally, from the definition of $\tau :$ $\pi _{1}(N,x_{0})\rightarrow
H_{1}(N),$ the map $\phi $ sends all elements $\overline{c_{i}}\in \tau
^{-1}([c_{i}])\subset \pi _{1}(N,x_{0})$ to $0\in \mathbb{Z},$ the elements $%
\overline{b_{i}}\in \tau ^{-1}([b_{i}])$ to $0\in \mathbb{Z}$ for $i\neq g$
and the elements $\overline{b_{g}}\in \tau ^{-1}([b_{g}])$ to $l\in \mathbb{Z%
}.$ So Claim 3 is proven.

Now we may construct a smooth submersion $f:N\rightarrow S^{1}$ such that,
if $f_{\ast }:\pi _{1}(N,x_{0})\rightarrow \pi _{1}(S^{1},y_{0})=\mathbb{Z}$
then $f_{\ast }=\phi =\psi \circ \tau .$ The construction of $f$ is
standard. We represent $N$ as a cell complex with one vertex $x_{0}.$ We
send this vertex to $y_{0}$ and map each cell of the 1-skeleton as
prescribed by $\phi .$ Notice here that the one-dimensional cells are
generators of the fundamental group. The map obtained can be extended to
two-dimensional cells which correspond to relations, and then to the
remaining three-dimensional cells, since $\pi _{2}(S^{1})=0.$

Assuming that $f$ is smooth and taking a regular value $y\in S^{1},$ $y\neq
y_{0}$ we consider the inverse image $f^{-1}(y).$ Then some component $T$ of 
$f^{-1}(y)$ is the required surface (see also Lemma 6.5 and 6.6 of \cite%
{Hempel} or Theorem 7.6 and Lemma 7.7 in \cite{Lackenbay}). In fact, under
isotopy, we may assume that $\partial T\cap (\cup _{i}c_{i})=\emptyset $
since $f$ sends each $c_{i}$ to a contractible curve in $S^{1}.$ Therefore
we may assume that the image $f(c_{i})$ avoid the value $y$ for each $i.$ On
the contrary, $f$ sends $b_{g}$ to a closed curve which wraps $l$ times
around $S^{1}.$ Additionally, the curve $b_{g}$ induces on a tubular
neighborhood of $T$ segments transverse to $T$ which are all transversely
oriented in a consistent way. This implies that $T$ does not separate $N.$
Therefore Proposition \ref{nice incompressible} is proven in Case I.

We consider now Case II.

In this case the kernel $Ker(\tau _{\#})$ of $\tau _{\#}:$ $H_{1}(\partial
N)\rightarrow H_{1}(N)$ is a non-empty subgroup of $H_{1}(\partial N)$ and
can be identified with a homogeneous linear system with unknown $[c_{i}]$
and $[b_{j}]$ and integer coefficients. These equations can be transformed
into equivalent ones by multiplying them by integers and adding them
together. Consequently, the linear system can be transformed into echelon
form. Below, $[c_{i}]$ and $[b_{j}]$ will also be referred to as variables
and we distinguish three subcases.

\textit{Subcase II(a): The last equation of the new system in echelon form
involves only the variables }$[c_{i}].$

In this subcase, there is an incompressible subcase $S$ in $N$ whose
boundary consists of curves which belong in the homotopy class of $c_{i}$ in 
$\partial N.$ We fill $N_{0}$ by gluing along each $S_{i}$ a solid torus $%
V_{i}$ such that $c_{i}$ is matched with a meridian of the solid torus $%
V_{i}.$ Thus, the properly embedded surface $S$ of $N$ is transformed to a
closed surface $S_{0}$ of $M.$

For $S$ we distinguish the following two possibilities:

(1) $S$ has genus $0.$

(2) $S$ has positive genus.

Thus, we first we assume that $S$ has genus $0.$ Then $S_{0}$ is a 2-sphere
and because of the irreducibility of $M,$ $S_{0}$ bounds a 3-ball $B$ in $M.$

If $\partial S$ has one boundary component then, without loss of generality,
we assume that $\partial S\subset S_{1}.$ We attach to each $S_{i},$ $%
i=2,\cdots ,n$ a solid torus $V_{i}$ such that each $c_{i}\in S_{i}$ is
matched with a meridian of $V_{i}.$ Let $M_{1}$ be the obtained manifold.
Obviously $\partial M_{1}=S_{1}.$ Then $S$ must be an essential disc in $%
M_{1}$ with $\partial S=c_{1}.$ But in this case $M_{1}$ should be a solid
torus and hence $M$ should be $S^{2}\times S^{1}$ which is impossible.

If $\partial S$ has has two boundary components then $S$ is a properly
embedded annulus in $N.$ But this is impossible since $N$ is acylindrical as
we have proved in Theorem \ref{attaching handlebody}.

If $\partial S$ has more than two boundary components then the existence of
the 3-ball $B$ implies that $S$ must have an even number of boundary
components. This follows from the fact that $S\subset N\subset N_{0}$ and
that each torus $S_{i}=\partial V_{i}$ which enters inside $B$ must exit of $%
B.$ Now, as $S$ has at least four boundary components we will show that $S$
is compressible in $N.$ In fact, let $c_{i}^{\prime },c_{i}^{\prime \prime }$
be two boundary components of $S$ belonging to some $S_{i}\cap B.$ Then each
one of the curves $c_{i}^{\prime },c_{i}^{\prime \prime }$ is parallel to $%
c_{i}$ within $S_{i}$ and so the curve of $S$ surrounding $c_{i}^{\prime }$
and $c_{i}^{\prime \prime }$ is trivial in $N\cap B.$

In conclusion, there is not an incompressible surface $S$ of genus $0,$
whose boundary components $c_{i}$ satisfy subcase II(a).

Second, we assume that $S$ has positive genus. Then, as in the previous
paragraph, we fill $N_{0}$ by gluing along each $S_{i}$ a solid torus $V_{i}$
such that $c_{i}$ is matched with a meridian of a solid torus $V_{i}.$ In
this way we get by construction the manifold $M.$ Again the properly
embedded surface $S$ of $N$ is transformed to a closed surface $S_{0}$ of $%
M. $ We compress $S_{0}$ in $M$ and either we obtain an incompressible
component, say $S_{0}^{\prime },$ in $M$ which is impossible by our
hypothesis or, we obtain a 2-sphere. In the latter case, we fall into the
scenario addressed earlier, which we have already proven to be impossible.
Therefore, Case II(\textit{a}) is entirely excluded.

\textit{Subcase II(b)}: The last equation of the new system is a equation
which relates more than one variable $[b_{j}]$ with the variables $[c_{i}].$

Note that in Subcase II(b) all classes $[c_{i}]$ may have coefficients equal
to $0,$ i.e. we have a relation only between the variables $[b_{j}].$

Thus, in Subcase II(b) we deduce that some $[b_{j}]$ is a free variable
independent of the variables $[c_{i}].$ Therefore, as in Claim 3, we may
construct a homomorphism $\phi $ which sends all $[c_{i}]$ to $0\in \mathbb{Z%
}$ and some $b_{j}$ to some integer $l\in \mathbb{Z}.$ Thus, we may
construct the required surface $T$ by considering again a smooth map $%
f:N\rightarrow S^{1}$ such that $f_{\ast }=\phi $ and choosing a component
of $f^{-1}(y)$ for a regular value $y$ of $f$ as in Claim 3.

\textit{Subcase II(c): The last equation of the new system is a relation
between a single variable }$[b_{j}],$\textit{\ say of }$[b_{1}],$\textit{\
and of }$[c_{i}].$

We remark that $[c_{1}]$ may or may not belong to the set of $[c_{i}]$ to
which $[b_{1}]$ is related via the preceding relation. In this situation, we
glue along each $S_{i},$ $i=2,\cdots ,n,$ a solid torus $V_{i}$ such that
each $c_{i}$ is matched with a meridian of $V_{i}.$ Let the resulting
manifold be denoted by $M_{1}.$ Obviously $\partial M_{1}=S_{1}.$ So the new
relation which relates $b_{1}$ with $c_{1}$ takes the form%
\begin{equation*}
k_{1}[b_{1}]+l_{1}[c_{1}]=0,~~k_{1},l_{1}\in \mathbb{Z},~~k_{1}\neq 0.
\end{equation*}

\textit{Claim 4}: The manifold $M_{1}$ is irreducible, atoroidal and
acylindrical.

\textit{Proof of Claim 4}. First, we observe that $M_{1}$ is irreducible.
This follows easily from the fact that $M$ is irreducible and the properties
of link $\mathcal{L}.$ Now, let $T$ be a torus in $M_{1}$ which is not
parallel to $S_{1}.$ We aim to show that that $T$ is compressible.
Specifically, we will assume that $T$ is incompressible in $M_{1}$ and
derive a contradiction. If $T\cap V_{i}\neq \emptyset $ for some $i\neq 1,$
then there will exist a simple closed curve $a\subset T\cap S_{i},$ where $%
S_{i}=\partial V_{i}.$ If $a$ is the boundary of a meridian of $V_{i}$ then,
from the irreducibility of $M_{1}$ we may assume that the curve $a$ is
non-trivial in $T.$ Therefore $a$ is an essential curve of $T$ which bounds
a disc in $M_{1}$ and hence $T$ should be compressible in $M_{1}$ which is
impossible since we have assumed that $T$ is incompressible in $M_{1}.$ We
know also that the only incompressible surfaces in $V_{i}$ are meridians of $%
V_{i}$ or annuli parallel to $S_{i}.$ Consequently, we deduce that $T\cap
(\cup _{i=2}^{n}S_{i})=\emptyset .$ On the hand, $T$ bounds a solid torus $V$
in $M$ and we claim that $V\subset M_{1}.$ In fact, if not then $%
V_{1}\subset V$ which is impossible from the properties of link $\mathcal{L}$
and the fact that $T\cap (\cup _{i=2}^{n}S_{i})=\emptyset .$ Therefore we
conclude that $T$ is compressible in $M_{1}.$

Finally, if $M_{1}$ were cylindrical then we could construct an
incompressible torus in $M_{1}$ which is non-parallel to $S_{1}.$ Thus,
Claim 4 is proved.

\textit{Claim 5}: We may find $m$ simple closed curves in $S_{1},$ all
parallel to a curve $d,$ such that:%
\begin{equation*}
m[d]=k_{1}[b_{1}]+l_{1}[c_{1}]=0.
\end{equation*}

\textit{Proof of Claim 5}. If $(k_{1},l_{1})=1,$ i.e. $k_{1}$ and $l_{1}$
are coprime, then there is a simple closed curve $d$ in $S_{1}$ that is
freely homotopic to the curve $b_{1}^{k_{1}}c_{1}^{l_{1}}$ of $S_{1}$ via a
homotopy entirely contained within $S_{1}.$ Otherwise, it will be $%
m(k_{1}^{\prime }[b_{1}]+l_{1}^{\prime }[c_{1}])=0,$\ where $(k_{1}^{\prime
},l_{1}^{\prime })=1.$\ In this latter case, there exists a simple closed
curve $d$\ in $S_{1}$\ such that $d$\ is freely homotopic to the curve $%
b_{1}^{k_{1}^{\prime }}c_{1}^{l_{1}^{\prime }}.$\ Hence we obtain $%
m[d]=m(k_{1}^{\prime }[b_{1}]+l_{1}^{\prime
}[c_{1}])=k_{1}[b_{1}]+l_{1}[c_{1}]=0.$\ This proves Claim 5.

Now there is a properly embedded surface $F$ in $M_{1}$ such that all its
boundary components are parallel copies of $d.$ Consider a parallel copy of $%
F,$ call it $F^{\prime }.$ Let $P=F\times \lbrack 0,1],$ where $F\times
\{0\} $ is identified with $F$ and $F\times \{1\}$ is identified with $%
F^{\prime }. $ Then $S_{1}\backslash P$ consists of $m$ annuli $A_{i}\subset
S_{1}.$ We form a closed surface $S$ in $M_{1}$ by gluing one boundary
component of $F$ with one boundary component of $F^{\prime }$ via an annulus
which is parallel to one of the $A_{i}.$ We claim that:

\textit{Claim 6}: By compressing $S$ in $M_{1},$ there are two possible
outcomes: either we obtain an incompressible surface $S_{0}$ in $M_{1}$ or $%
M_{1}$ is homeomorphic to a surface bundle over $S^{1}.$

\textit{Proof of Claim 6}. If we do not arrive at an incompressible $S_{0},$
then we deduce that $Q=M_{1}\backslash P$ is a handlebody, which is
necessarily homeomorphic to a trivial product $F\times \lbrack 0,1].$ Since
the manifolds $P$ and $Q$ are glued along parts of $\partial P$ and $%
\partial Q$ which are homeomorphic to $F,$ we deduce that $M_{1}$ is
homeomorphic to a surface bundle over $S^{1}$ with fiber homeomorphic to $F$
and thus Claim 6 is proved.

Now since $M_{1}$ has one toral component we deduce that $F$ has one
boundary component in the case that $M_{1}$ is a surface bundle over $S^{1}.$
Thus, we examine each of the following two cases separately:

(1) $S_{0}$ is incompressible in $M_{1}.$

Then, since $S$ in neither an essential disc nor an annulus in $M_{1},$ as $%
M_{1}$ is acylindrical, we deduce that $S_{0}$ has positive genus greater
than 1 and we will show that $S_{0}$ is incompressible in $M.$ In fact, if
not, then there is an essential curve $d$ of $S_{0}$ which bounds a disc $D$
in $M$ and $D\cap S_{0}=d.$ Since $S_{0}$ is incompressible in $M_{1}$ it
must be $D\cap S_{1}\neq \emptyset .$ Therefore, there must be a curve $%
d_{0} $ in $D\cap S_{1}$ which bounds a disc $D_{0}$ such that either $%
D_{0}\subset V_{1}$ or $D_{0}\subset M_{1}$ which is impossible. Note here
that only the curve $c_{1}$ on $S_{1}$ is contractible in $V_{1}$ and $d_{0}$
cannot be isotopic to $c_{1}.$

(2) $M_{1}$ is homeomorphic to a surface bundle over $S^{1}.$

In this case we will prove that $M$ has a finite covering which is Haken
i.e. $M$ is virtual Haken. First we notice that $M_{1}=(F\times I)/\phi ,$
where $\phi $ is pseudo-Anosov. In fact, if $\phi $ was periodic or
reducible we may prove that $M_{1}$ is toroidal contradicting in this way
Claim 4.

Now, we have the following claim.

\textit{Claim 7}: We may find a finite covering $\widetilde{M_{1}}$ of $%
M_{1} $ such that the curve $d$ lifts to $l_{1}$ copies, say $\widetilde{d},$
such that $\widetilde{d}=\widetilde{c_{1}}+k_{1}\widetilde{b_{1}},$ where $%
\widetilde{c_{1}},$ $\widetilde{b_{1}}$ are lifts of $c_{1},$ $b_{1}$ in $%
\widetilde{M_{1}}.$

\textit{Proof of Claim 7}. First we may find a finite cover $\widetilde{S_{1}%
}$ of $S_{1}=\partial M_{1}$ such that: the curve $c_{1}$ lifts to a curve $%
\widetilde{c_{1}}$ which covers $l_{1}$ times the curve $c_{1};$ the curve $%
b_{1}$ lifts to $l_{1}$ copies, say $\widetilde{b_{1}}\subset \widetilde{%
S_{1}}.$ So, the curve $d$ lifts to $l_{1}$ copies of a curve $\widetilde{d}%
\subset \widetilde{S_{1}}$ such that $\widetilde{d}=\widetilde{c_{1}}+k_{1}%
\widetilde{b_{1}}.$

The fibers $F$ of $M_{1}=(F\times I)/\phi $ induce on $S_{1}=\partial M_{1}$
a trivial product by parallel copies of $d.$ Now we may construct a finite
covering $\widetilde{M_{1}}$ of $M$ by extending the covering map $q:$ $%
\widetilde{S_{1}}\rightarrow S_{1}$ to a covering map $q^{\prime }:$ $%
\widetilde{M_{1}}\rightarrow M_{1}$ as follows: $\widetilde{M_{1}}$ will be
defined again as a fibration of the form $\widetilde{M_{1}}=(\widetilde{F}%
\times I)/\phi $ such that the fibers $\widetilde{F}$ are homeomorphic to $F$
and they induce on $\widetilde{S_{1}}=\partial \widetilde{M_{1}}$ a trivial
product by curves $\widetilde{d}_{t}$ parallel to $\widetilde{d}.$ In fact,
it is sufficient to consider the product $\widetilde{F}\times I$ such that $%
\partial (\widetilde{F}\times \{t\})=\widetilde{d}_{t},$ $t\in I,$ and
define $\widetilde{M_{1}}=(\widetilde{F}\times I)/\phi .$ Therefore the
covering map $q:$ $\widetilde{S_{1}}\rightarrow S_{1}$ can be extended to a
covering map $q^{\prime }:$ $\widetilde{M_{1}}\rightarrow M_{1}$ which has
all required properties. This proves Claim 7.

From Claim 7, we glue a solid torus $V$ along $\partial \widetilde{M_{1}}$
so that $\widetilde{d}$ is matched to a meridian of $V.$ Call the resulting
manifold $\widetilde{M}.$ Then $\widetilde{M}$ is homeomorphic, via some
homeomorphism $h,$ to a manifold $M^{\prime }$ obtained by gluing a solid
torus $V$ along $\partial \widetilde{M_{1}}$ in such a way that $\widetilde{%
c_{1}}$ is matched to a meridian of $V.$ Therefore $\widetilde{M},$ and
hence $M^{\prime },$ is a finite covering of $M.$ Moreover, in $\widetilde{M}
$ we may find an incompressible surface $\widetilde{F}.$ This is actually
the surface $F$ in which a disc $D$ is glued along $\partial F.$ Therefore $%
h(\widetilde{F})$ is an incompressible surface in $M^{\prime }.$
\end{proof}

The following remark is important for the remainder of this paper.

\begin{remark}
\label{separating boundary}The properly embedded incompressible surface $T$
constructed in Proposition \ref{nice incompressible} is, by construction,
non-separating in $N$ but its boundary $\partial T$ may separate the
boundary $\partial N$ of $N.$
\end{remark}

\begin{center}
\begin{figure}[h]
%\hspace*{9mm}
\includegraphics[scale=0.85]
%, bb = 46 316 565 633]
{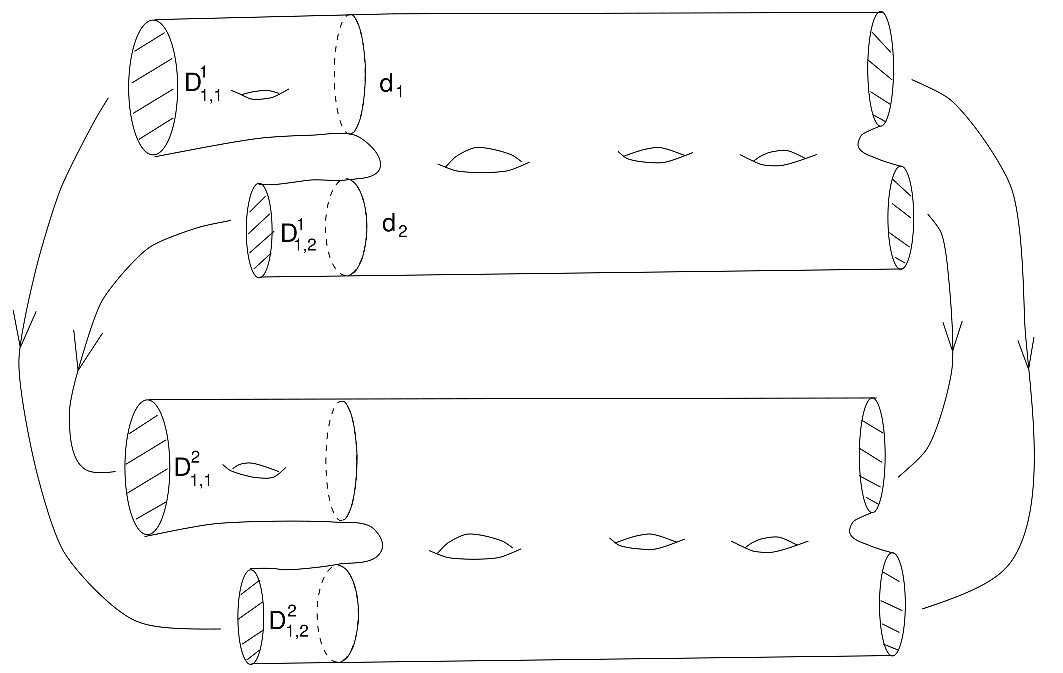}
%\includegraphics[scale=0.7]{HAKEN_9.pdf} %\\vspace*{-8mm}\newline
%\label{disk4}
\caption{The double-covering of a handlebody.}
\end{figure}
\end{center}

\section{Construction of the appropriate covering space}

In this section we will adapt ideas from the proof of Theorem 1.1 of \cite%
{Cooper-Long-Reid}, combined with Proposition \ref{nice incompressible}, to
establish the framework for proving the virtual Haken conjecture.
Specifically, we will construct a finite cover of $M$ which contains a
separating incompressible surface.

We note that by Corollary \ref{WLG}, $N$ admits a complete hyperbolic
structure as in Theorem 1.1 of \cite{Cooper-Long-Reid}, though this fact is
not essential to the subsequent proof.

We recall that, up to now, we have constructed a compact, orientable,
incompressible surface $T,$ properly embedded in $N.$ By construction the
surface $T$ is non-separating in $N;$ however, we cannot control its
boundary, so $\partial T$ may separate the boundary $\partial N$ of $N.$ In
that case, our goal will be to eliminate this phenomenon by passing to a
finite (in fact, double) covering of $N$ and of the closed manifold $M;$
recall that $M$ is obtained from $N$ by attaching a handlebody $H$ along $%
\partial N$ (of genus $n$, say), as per our previous considerations.

To clarify the construction of the double covering $N^{\prime }$ of $N,$ we
first construct a double covering $Q^{\prime }$ of a handlebody $Q$ of index 
$n\geq 2.$ This preliminary construction illustrates the basic idea
underlying the construction of $N^{\prime }.$

Let $Q$ be a handlebody of genus $n\geq 2,$ and let $d$ be a multicurve of $%
\partial Q$ (possibly connected) that separates $\partial Q$ into two
components. Each component $d_{i},$ $i=1,...,m$ of $d$ is assumed to be a
meridional curve, that is, the boundary of a meridian disk of $Q.$ Let $%
Q_{1},$ $Q_{2}$ be the two components separated by $d.$ For $i=1,2,$ we
choose a meridional curve $a_{i}$ in $Q_{i}$ that does not separate $Q_{i}.$
Let $D_{i}$ be a meridian of $Q$ with $\partial D_{i}=a_{i}.$ We cut and
open $Q$ along $a_{1}$ and $a_{2},$ and we denote by $Q_{0}$ the resulting
space, which is a "spotted" handlebody. Thus, the curve $a_{1}$ (resp. $%
a_{2})$ gives rise to two curves $a_{1,1},$ $a_{1,2}$ (resp. $a_{2,1},$ $%
a_{2,2})$ in $Q_{0}.$ We then consider two copies $Q_{0}^{1},$ $Q_{0}^{2}$
of $Q_{0}$ and for $i,j\in \{1,2\},$ denote by $a_{i,j}^{1}$ (resp. $%
a_{i,j}^{2})$ the curves of $Q_{0}^{1}$ (resp. $Q_{0}^{2})$ corresponding to 
$a_{i,j}$ of $Q_{0},$ see Figure 9. By $D_{i,j}^{k}$ we denote the meridians
bounded by $a_{i,j}^{k},$ where $i,j,k\in \{1,2\}.$ Now we construct a
double covering $Q^{\prime }$ of $Q$ by matching together the discs

\begin{itemize}
\item $D_{1,1}^{1}$ with $D_{1,2}^{2},$ $D_{1,2}^{1}$ with $D_{1,1}^{2},$ $%
D_{2,1}^{1}$ with $D_{2,2}^{2},$ $D_{2,2}^{1}$ with $D_{2,1}^{2}$, see
Figure 9, where the arrows show the meridians which are matched together.
\end{itemize}

In this way, we obtain a handlebody $Q^{\prime }$ of genus $n^{\prime }>n,$
which is a double covering space of $Q.$ This construction will be referred
to as the \textit{double covering} of $Q.\smallskip $

By Proposition \ref{nice incompressible} there exists a properly embedded,
non-separating surface $T$ in $N$ such that $\partial T\cap \mathcal{S}%
=\emptyset ,$ where $\mathcal{S}=\{s_{i},i=1,...,n\}.$ The curves $s_{i}$
are meridional in $H.$ Therefore, each component of $\partial T$ is also a
meridional curve in $H.$ Indeed, each $s_{i}$ is a meridional curve of $H$
and by cutting $H$ along the meridians bounded by $s_{i},$ one obtains a
3-ball $B.$ The boundary components of $\partial T$ then lie on the boundary
of this ball $B,$ which forces them to bound meridians in $H.$ As we noted
in Remark \ref{separating boundary}, the boundary $\partial T$ may be a
curve or a multicurve that separates $\partial N=\partial H.$ In this case,
we will find a double covering $p:M^{\prime }\rightarrow M$ such that:

\begin{itemize}
\item[$(I)$] If $N^{\prime }=p^{-1}(N)$ and $H^{\prime }=p^{-1}(N)$ then both
$N^{\prime }$ and $H^{\prime }$ are connected;

\item[$(II)$] we can choose a component $T^{\prime }$ of $p^{-1}(T)$ which
does not separate $M^{\prime }$ and whose the boundary $\partial T^{\prime }$
also does not separate $\partial M^{\prime }.$
\end{itemize}

Therefore our next goal is to prove the following proposition.

\begin{proposition}
\label{good non-separating cover}We assume that $\partial T$ separates $%
\partial N=\partial H.$ Then, we may pass to a double covering $p:M^{\prime
}\rightarrow $ $M$ so that the preceding conditions $(I)$ and $
(II)$ hold.
\end{proposition}

\begin{proof}
The properly embedded, connected, incompressible surface $T$ in $N$ gives
rise to a closed surface $\overline{T}$ is $M.$ This follows from the fact
that each component of $\partial T$ is a meridional curve in $H,$ and hence $%
\overline{T}$ is formed by adding to $T$ meridians of $H$ along each
component of $\partial T.$ The surface $T$ is non-separating in $N;$
therefore $\overline{T}$ is non-separating in $M.$ By cutting open $M$ along 
$\overline{T},$ we obtain a connected manifold $\overline{M}$ with two
boundary components $\overline{T}_{1},$ $\overline{T}_{2},$ each
corresponding to $\overline{T}.$

Since $\partial T$ separates $\partial H,$ in the interior of $\overline{M}$
there are two distinct "spotted" handlebodies $H_{1},$ $H_{2}$ arising from
splitting $H$ along the meridians bounded by the components of $\partial T.$
For each $i=1,2,$ the spots of $H_{i}$ arise exactly from these cutting
meridians. Clearly, the spots of $H_{i}$ lie in the respective boundary
component $\overline{T}_{i}.$

Choose a non-separating meridian $D_{1}$ in $H_{1}\cap Int(\overline{M})$
and a non-separating meridian $D_{2}$ in $H_{2}\cap Int(\overline{M}).$ Let $%
a_{i}=\partial D_{i}$ for $i=1,2.$ Let%
\begin{equation*}
\overline{N}=\overline{M}\backslash (Int(H_{1})\cup Int(H_{2})).
\end{equation*}%
For each $i=1,2,$ considering an annular neighborhood $U(a_{i})\subset
\partial \overline{N}$ of $a_{i}$ denote by $a_{i,1}$ and $a_{i,2}$ the
boundary components of $U(a_{i}).$ In $\overline{N}$ we may construct two
distinct properly embedded, surfaces $F_{1}$ and $F_{2}$ such that: $%
\partial F_{i}=a_{i,1}\cup a_{i,2}$ for $i=1,2.$ Note that each $F_{i}$ need
not be connected. More precisely, the surface $F_{i}$ is given by the
following union 
\begin{equation*}
F_{i}=(\partial H_{i}\backslash U(a_{i}))\cup (T_{i}\backslash (\text{spots
of }H_{i}\text{ in }T_{i})),~i=1,2,
\end{equation*}%
and we may then isotope each $F_{i}$ so that it is properly embedded in $%
\overline{N}$ with $\partial F_{i}=a_{i,1}\cup a_{i,2}.$

In the manifold $M,$ obtained by gluing back $T_{1}$ with $T_{2},$ we attach
to $F_{i}$ parallel meridians $D_{i,1},$ $D_{i,2}$ of $H_{i},$ where $%
\partial D_{i,j}=a_{i,j},$ $i,j\in \{1,2\},$ and denote by $\overline{F}_{i} 
$ the resulting surface corresponding to $F_{i}$ (for $i=1,2).$ Thus $%
\overline{F}_{1},\overline{F}_{2}\subset M$ and by construction, they are
distinct, and each one is a non-separating surfaces in $M.$ The latter
follows from the fact that $T$ is non-separating in $N.$

We will proceed now to the construction of the double covering space $%
M^{\prime }\rightarrow M.$ This construction is carried out in parallel with
the double-covering of a handlebody described earlier. First, we cut open $M$
along $\overline{F}_{1}$ and $\overline{F}_{2}.$ Doing so yields a manifold $%
M_{0}$ with four boundary components: two coming from $\overline{F}_{1},$
which we label $\overline{F}_{1,1}$ and $\overline{F}_{1,2},$ and two from $%
\overline{F}_{2}$ which we label $\overline{F}_{2,1}$ and $\overline{F}%
_{2,2}.$ Consider now two copies of $M_{0}^{1}$ and $M_{0}^{2}$ of $M_{0}.$
The components corresponding to $\overline{F}_{i,j},$ for $i,j\in \{1,2\},$
which belong to $M_{0}^{1}$ (resp. $M_{0}^{2})$ are labeled by $\overline{F}%
_{i,j}^{1}$ (resp. $\overline{F}_{i,j}^{2}).$ For $j=1,2,$ the meridians $%
D_{1,j}$ $\subset \overline{F}_{1}$ and $D_{2,j}\subset \overline{F}_{2}$
defined previously, induce discs $D_{i,j}^{k}$ in $\overline{F}_{i,j}^{k}$
for each $i,j,k\in \{1,2\},$ in accordance with the terminology introduced
above in the construction of the double-covering of a handlebody.

Finally, we glue $M_{0}^{1}$ with $M_{0}^{2}$ along their boundaries so that:

\begin{itemize}
\item the surface $\overline{F}_{1,1}^{1}$ is glued with the surface $%
\overline{F}_{1,2}^{2}$ in such a way that $D_{1,1}^{1}$ is matched with the
surface $D_{1,2}^{2};$

\item the surface $\overline{F}_{1,2}^{1}$ is glued with the surface $%
\overline{F}_{1,1}^{2}$ in such a way that $D_{1,2}^{1}$ is matched with the
surface $D_{1,1}^{2};$

\item the surface $\overline{F}_{2,1}^{1}$ is glued with the surface $%
\overline{F}_{2,2}^{2}$ in such a way that $\overline{D}_{2,1}^{1}$ is
matched with the surface $D_{2,2}^{2};$

\item the surface $\overline{F}_{2,2}^{1}$ is glued with the surface $%
\overline{F}_{2,1}^{2}$ in such a way that $D_{2,2}^{1}$ is matched with the
surface $D_{2,1}^{2}.$
\end{itemize}

We then obtain a closed manifold $M^{\prime },$ which is a double-covering
of $M.$ Furthermore, in the interior of $M^{\prime }$ there is a handlebody $%
H^{\prime }$ which is a double cover of $H$ because it is built following
the same pattern as the double-covering construction for a handlebody.
Clearly, $N^{\prime }=M^{\prime }\backslash H^{\prime }$ is a double
covering of $N=M\backslash H.$ If $p:M^{\prime }\rightarrow M$ denotes the
aforementioned double covering, then by construction we may choose a
component $T^{\prime }$ of $p^{-1}(T)$ which does not separate $M^{\prime }$
and and its boundary $\partial T^{\prime }$ likewise does not separate $%
\partial M^{\prime }.$ This finish the proof of the proposition.
\end{proof}

From Proposition \ref{good non-separating cover}, we may, after possibly
passing to a finite covering of $M,$ assume that%
\begin{equation}
\partial T\text{ does not separate }\partial M.\smallskip  \tag{$A$}
\end{equation}%
Going forward, we adopt this assumption in all subsequent arguments.

\begin{proposition}
\label{essential3}Assume that $M$ is not Haken and that $M$ has no finite
cover which is a surface bundle over $S^{1}.$ Then the manifold $%
N=M\backslash (Int(H))$ contains an essential closed surface $S$ such that:

(1) $S$ can be lifted to an embedded, incompressible, separating surface $%
\widetilde{S}$ in a finite cover $p:$ $\widetilde{N}\rightarrow $ $N$ of
order $d.$

(2) The preimage under $p$ of each curve $s_{i}\subset $ $\widetilde{N}$
consists of $d$ closed curves $\widetilde{s}_{i}^{1},...,\widetilde{s}%
_{i}^{d};$ moreover, for each $i,j,$ the restriction 
\begin{equation*}
p_{|\widetilde{s}_{i}^{j}}:\widetilde{s}_{i}^{j}\rightarrow s_{i}
\end{equation*}%
is a homeomorphism.
\end{proposition}

\begin{proof}
By Proposition \ref{nice incompressible} there exists a properly embedded,
non-separating surface $T$ in $N$ such that $\partial T\cap \mathcal{S}%
=\emptyset .$ Also, according to Assumption $(A),$ we may assume that $%
\partial T$ does not separate $\partial M.$

Let $U(T)$ denote a tubular neighborhood of $T,$ and define $X=N\backslash
Int(U(T)).$ Let $X_{i}=X\times \{i\},$ $i\in \mathbb{Z}$ be a copy of $X.$
We consider the infinite cyclic cover%
\begin{equation*}
p:N_{T}\rightarrow N
\end{equation*}%
dual to $T,$ where $N_{T}=\cup _{i}X_{i}$ and $X_{i}\cap X_{i+1}=T_{i}$ is a
copy of $T.$ Since $T$ is incompressible in $N$ it follows that $T_{i}$ is
incompressible in $N_{T}.$ Define the submanifold $Y_{k}=\cup
_{i=1}^{k}X_{i} $ of $N_{T}$ and obviously we have that 
\begin{equation*}
N_{T}=E_{-}\cup _{T_{0}}Y_{k}\cup _{T_{k}}E_{+}
\end{equation*}%
where $E_{-}=\cup _{i\leq 0}X_{i},$ $E_{+}=\cup _{i\geq k+1}X_{i}.$

Let $W=\partial Y_{k}\backslash (T_{0}\cup T_{k+1})$ and we consider the
surface $S^{\prime }=W^{\prime }\cup T_{0}^{\prime }\cup T_{k+1}^{\prime },$%
\ where $W^{\prime }$\ is parallel to $W\subset \partial Y_{k}$\ and $%
T_{0}^{\prime },$\ $T_{k+1}^{\prime }$\ are parallel to $T_{0},$\ $T_{k+1}$\
respectively in $N_{T}.$

The finite covering $\widetilde{N}$ is defined as the quotient manifold $%
Y_{k}/\sim ,$ where the two copies $T_{0}$ and $T_{k+1}$ of $T$ in $Y_{k}$
are identified via the identity map. We note that we may choose $k$
arbitrarily large in the definition of $Y_{k}.$ Now, by compressing $%
S^{\prime }$ in $\widetilde{N}$ \textbf{\ }we obtain a surface $\widetilde{S}
$ in $\widetilde{N}.$ Clearly, by construction $S^{\prime }$ is separating
in $\widetilde{N},$ and hence the same is true for $\widetilde{S}.$
Projecting $S$ into $N$ via the covering projection $\widetilde{p}:%
\widetilde{N}$ $\rightarrow N$ we take the required surface $S$ in $N,$ i.e. 
$S=\widetilde{p}(\widetilde{S}).$ First, we observe that the boundaries of
the compressing discs, along which $S^{\prime }$ is compressed to $%
\widetilde{S},$ consist of arcs lying in $W^{\prime }$ as well as in $%
T_{0}^{\prime }$ and $T_{k+1}^{\prime },$ since both $W$ and $T$ are
incompressible in $N.$ Second, we note that since $n$ can be chosen
arbitrarily large, the number of disjoint properly embedded arcs in $%
T_{0}^{\prime }\cup T_{k+1}^{\prime }$ is sufficiently smaller than the
number of disjoint properly embedded arcs in $W.$ Consequently, the
compression process cuts $W^{\prime }$ into components of arbitrarily large
genus. This implies that we can obtain a closed connected surface $S$ in $%
\widetilde{N}$ of genus at least 2, which is incompressible in $\widetilde{N}%
.$

The final step is to prove that $\widetilde{S}$\ is not parallel to $%
\partial \widetilde{N}.$ Assuming the contrary, we would reach a
contradiction as follows. First, we may assume that the boundaries of the
compressing discs along which $S^{\prime }$ is compressed to $\widetilde{S}$
do not intersect each other. Therefore, the boundary $d$\ of a cutting disc $%
D$ of \ $S^{\prime }$ contains arcs $a_{j}$ that are properly embedded
essential arcs of $T_{0}^{\prime }\cup T_{k+1}^{\prime }.$ \ Recall that, by
construction, the components of $\partial (T_{0}^{\prime }\cup
T_{k+1}^{\prime })$ are curves of $\partial \widetilde{N}$ (or more
precisely parallel to curves of $\partial \widetilde{N}).$ Thus, if $%
\widetilde{S}$\ were parallel to $\partial \widetilde{N},$ then,in a trivial
product neighborhood of $\partial \widetilde{N},$ the curve $d$\ would be
parallel to the boundary of a disc $E\subset \partial \widetilde{N},$ since $%
\partial \widetilde{N}$ is incompressible in $\widetilde{N}.$ This \ implies
that each arc $a_{j}\subset T_{0}^{\prime }\cup T_{k+1}^{\prime }$ would not
be a properly embedded essential arc of $T_{0}^{\prime }\cup T_{k+1}^{\prime
},$ contradicting our assumption.
\end{proof}

The fact that $\partial T$ does not separate $N$ is crucial because it
ensures that $\partial \widetilde{N}$ is connected. Moreover, by
construction the covering $\widetilde{p}:\widetilde{N}$ $\rightarrow N$ is
cyclic and this implies that the induced covering%
\begin{equation*}
\widetilde{p}_{|\partial \widetilde{N}}:\partial \widetilde{N}\rightarrow
\partial N
\end{equation*}%
is also cyclic.

Let $n$ be the genus of $\partial N,$ $b$ the number of boundary components
of $T,$ and $k$ the degree of the cyclic covering $\widetilde{p}:\widetilde{N%
}$ $\rightarrow N.$ Below we will compute the genus $g$ of $\partial 
\widetilde{N}$ by viewing $\partial \widetilde{N}$ as a cyclic covering of $%
\partial N.$ In fact, $\partial \widetilde{N}$ is built from  $k$ copies $%
X_{i}$ of a surface of genus  $n-b,$ each copy having $2b$ boundary
components. These copies are glued cyclically: $b$ boundary components of $X_{i}$ are glued to $b$ boundary components of $X_{i+1({mod}\,\,k)}.$ One
then checks that 
\begin{equation*}
g=k(n-b)+k(b-1)+1=kn-k+1=k(n-1)+1.
\end{equation*}

Since $\partial T\cap \mathcal{S}=\emptyset ,$ where $\mathcal{S}=\{s_{i},$ $%
i=1,\cdots ,n\},$ we deduce, as we have already noted, that each boundary
component of $\partial T$ bounds a meridian in the handlebody $H$ that is
attached to $N$ to form $M.$ Let us denote by $\Delta $ the union of
meridians in $H$ bounded by the boundary components of $T.$

Consequently, the same pattern of assembling the copies $X_{i}$ and which
describes the construction of cyclic coverings $\widetilde{p}:\widetilde{N}$ 
$\rightarrow N$ and $\widetilde{p}_{|\partial \widetilde{N}}:\partial 
\widetilde{N}\rightarrow \partial N$ can be carried over to $H.$ Indeed, if
we denote by $\Delta _{i}$ the union of discs bounded by the boundary
components of $\partial T_{i}$ (which corresponds to $\Delta ),$ each
surface $X_{i}$ can be replaced by a "spotted" handlebody $H_{_{i}},$ whose
"spots" are exactly the union of sets $\Delta _{i-1}$ and $\Delta _{i}.$
Therefore, we have 
\begin{equation*}
\partial H_{i}=(\partial X_{i}\backslash (T_{i-1}\cup T_{i}))\cup \Delta
_{i-1}\cup \Delta _{i}.
\end{equation*}%
Gluing these "spotted" handlebodies in the same cyclic fashion as with the $%
X_{i},$ one obtains a cyclic covering. 
\begin{equation*}
p_{H}:\widetilde{H}\rightarrow H.
\end{equation*}%
Clearly, the restriction to the boundary,%
\begin{equation*}
p_{H|\partial \widetilde{H}}:\partial \widetilde{H}\rightarrow \partial H
\end{equation*}%
coincides exactly with $\widetilde{p}_{|\partial \widetilde{N}}:\partial 
\widetilde{N}\rightarrow \partial N.$

Finally, the above analysis implies that the two coverings $\widetilde{p}:%
\widetilde{N}$ $\rightarrow N$ and $p_{H}:\widetilde{H}\rightarrow H$ match
together, so that we obtain the following corollary.

\begin{corollary}
\label{Indeed finite cover} There exists a finite covering $\widetilde{M}%
\rightarrow M$ such that $\widetilde{M}=$ $\widetilde{N}$ $\cup _{\partial 
\widetilde{N}=\partial \widetilde{H}}\widetilde{H}.$
\end{corollary}

\begin{proof}
The finite covering $\widetilde{M}\rightarrow M$ arises by matching the
coverings $\widetilde{p}:\widetilde{N}$ $\rightarrow N$ and $p_{H}:%
\widetilde{H}\rightarrow H$ along their common boundaries.
\end{proof}

\begin{remark}
The covering $\widetilde{M}\rightarrow M$ obtained above is not cyclic in
the case that $\partial T$ does not separate $\partial N.$
\end{remark}

\section{Proof of the virtual Haken conjecture}

Now we are able to proceed with the proof of the Virtual Haken Conjecture,
which can be formulated in the following theorem.

\begin{theorem}
\label{Haken conjecture final} Every compact, orientable, irreducible
three-dimensional manifold $M$ with infinite fundamental group is virtually
Haken. That is, it has a finite cover which is a Haken manifold.
\end{theorem}

\begin{center}
\begin{figure}[h]
%\hspace*{9mm}
\includegraphics[scale=0.8, bb = 46 316 565 633]{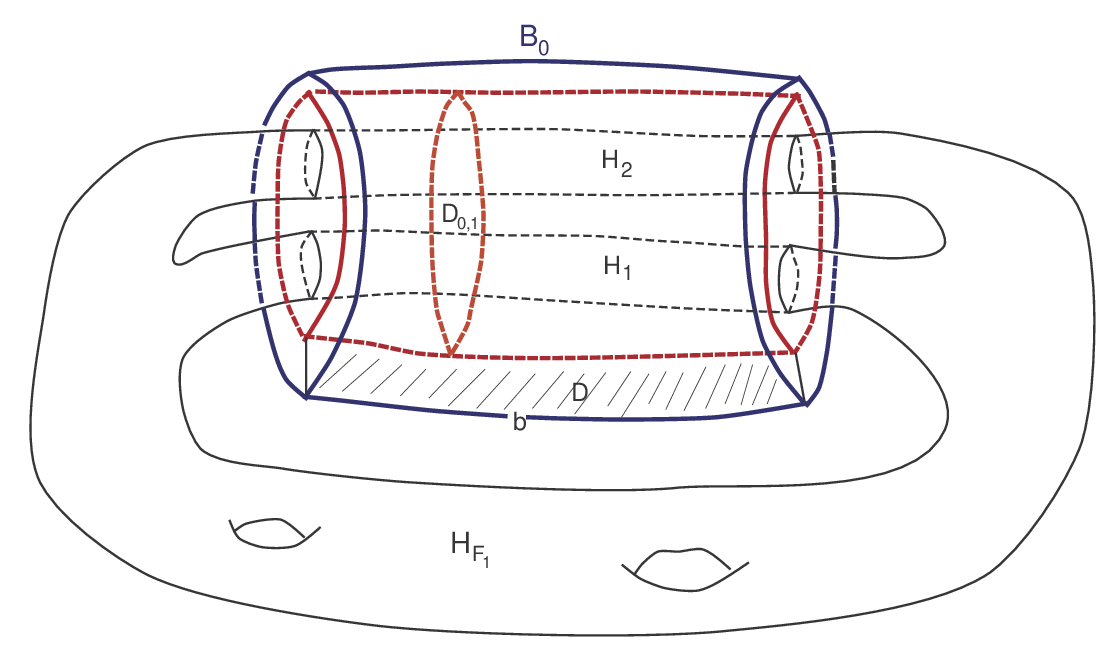}
%\includegraphics[scale=0.7]{HAKEN_9.pdf} %\\vspace*{-8mm}\newline
%\label{disk4}
\caption{A curve $b$ of $\widetilde{S}$ which bounds a disc $D\subset 
\widetilde{N}$ in Subcase 2(a).}
\end{figure}
\end{center}

\begin{proof}
\noindent We assume that $M$ is non-Haken and that $M$ does not have a
finite covering which is a surface bundle over $S^{1}.$ By Theorem \ref%
{attaching handlebody} there is a 3-manifold $N$ with connected
incompressible boundary such that $M$ is obtained by filling $N$ by a
handlebody of genus $\geq 2.$ Then, by Proposition \ref{essential3} and
Corollary \ref{Indeed finite cover}, there exists a finite covering $%
\widetilde{M}\rightarrow M$ in which lies a cyclic finite covering $%
\widetilde{N}\rightarrow N$ of $N;$ furthermore in $\widetilde{N}$ there is
a closed, separating surface $\widetilde{S}$ which is not parallel to $%
\partial \widetilde{N}.$

If $\widetilde{S}$ is not incompressible in $\widetilde{M},$ we continue to
compress $\widetilde{S}$ in $\widetilde{M}$ until we obtain a surface $%
\widetilde{S}_{0}.$ We will show that one component of $\widetilde{S}_{0}$
is incompressible and of genus $g\geq 2$ in $\widetilde{M}.$

Recall that $\widetilde{M}$ is obtained from $\widetilde{N}$ by gluing a
handlebody $\widetilde{H}$ to $\widetilde{N}$ via a homeomorphism $%
h:\partial \widetilde{H}\rightarrow \partial \widetilde{N}$ and let $%
\widetilde{\mathcal{L}}$ denote the lifting of link $\mathcal{L}$ in $%
\widetilde{M}$ which is obviously contained in $\widetilde{H}.$ Since $%
\widetilde{S}$ is assumed compressible in $\widetilde{M},$ choose a
compressing disc $D\subset \widetilde{M}$ such that $D\cap \widetilde{S}%
=\partial D.$ Because $\partial \widetilde{N}$ is incompressible in $%
\widetilde{N},$ $D$ decomposes as a planar surface $P\subset \widetilde{N}$
together with sub-discs $D_{i}\subset \widetilde{H}.$ By compressing $%
\widetilde{S}$ successively in $\widetilde{M},$ we eventually obtain a
surface (possibly disconnected) $\widetilde{S}_{0},$ and simultaneously $%
\partial \widetilde{H}$ is cut along the discs $D_{i}$ to a surface $F$
(possibly disconnected). Denote by $F_{1},\cdots ,F_{m}$ the components of $%
F;$ for each $i,$ $F_{i}$ is the boundary of a handlebody $H_{F_{i}}.$
Clearly, to obtain the components $H_{F_{i}},$ one removes from $\widetilde{H%
}$ a collection of 3-balls $B_{1},...,B_{k},$ each $B_{i}$ being
homeomorphic to $D_{i}\times \lbrack 0,1].$

We distinguish the following two cases:

(1) $\widetilde{S}_{0}$ is not simply connected. In that case, select one of
its non-simply-connected components (which we will continue to denote by $%
\widetilde{S}_{0}),$ which is incompressible in $\widetilde{M}.$

(2) Each component of $\widetilde{S}_{0}$ is a sphere. Then we may assume
that $\widetilde{S}_{0}$ is connected. Suppose for contradiction that it is
not. During the compression process of $\widetilde{S},$ let $X$ be a
compressing disc, whose compression along it, yields two copies $X_{1},$ $%
X_{2},$ lying in two distinct components $\widetilde{S}_{0}^{1},$ $%
\widetilde{S}_{0}^{2}$ of $\widetilde{S}_{0}.$ Since each $\widetilde{S}%
_{0}^{1},$ $\widetilde{S}_{0}^{2}$ is a 2-sphere we can glue them back
together along $X_{1},$ $X_{2}$ to form a single 2-sphere. Thus we reduce the
number of components of $\widetilde{S}_{0}.$ Therefore, we may assume that
there do not exist such discs $X_{1},$ $X_{2}$ in different components $%
\widetilde{S}_{0}^{1},$ $\widetilde{S}_{0}^{2}.$ In particular, any
separating compressing disc $X$ of $\widetilde{S}$ may not be used in the
compression process. Therefore, the compression process can be arranged so
that $\widetilde{S}_{0}$ remains connected, establishing our claim.

Now, since $\widetilde{M}$ is irreducible, $\widetilde{S}_{0}$ bounds a
3-ball $B_{0}$ in $\widetilde{M}.$ The components $F_{1},\cdots ,F_{m}$ of $%
F $ lie either in $B_{0}$ or in its complement $B_{0}^{c}=\widetilde{M}%
\backslash B_{0}.$ This follows from the fact that both $\widetilde{S}$ and $%
\partial \widetilde{N}$ are connected. We distinguish two subcases, which
will examined separately.

(2a) The components $H_{F_{i}}$ are contained in $B_{0}^{c}.$

Then we can show that $\widetilde{S}$ must be compressible in $\widetilde{N}%
, $ contradicting our assumption. Indeed, $B_{0}$ can be viewed as a product 
$D^{2}\times \lbrack 0,1],$ i.e. $B_{0}=$ $D^{2}\times \lbrack 0,1].$ The
handlebody $\widetilde{H}$ and the surface $\widetilde{S}$ can be recovered
as follows:

In $B_{0}$ we construct 3-balls $H_{1},\cdots ,H_{k},$ each $H_{i}$ being $%
D_{i}\times \lbrack 0,1]\subset D^{2}\times \lbrack 0,1]$ for some 2-disc $%
D_{i}\subset D^{2},$ with $D_{i}\times \{0,1\} \subset D^{2}\times \{0,1\}.$
Now in $B_{0}$ we \textquotedblleft open holes\textquotedblright \ by
removing the 3-balls $H_{i}.$ Clearly $B_{0}\backslash (\cup _{i}H_{i})$ is
a handlebody. The handlebody $\widetilde{H}$ is recovered by attaching the
3-balls $H_{i}\subset B_{0}$ to the handlebodies $H_{F_{j}}\subset B_{0}^{c}$
along the discs $D_{i}\times \{0,1\} \subset H_{i}.$ In Figure 10, a
handlebody $H_{F_{1}}$ in $B_{0}^{c}$ and balls $H_{1},H_{2}$ in $B_{0}$ are
illustrated.

To recover the surface $\widetilde{S}$ we consider 2-discs $D_{0,n},$ $%
n=1,...,a$ such that:

\begin{itemize}
\item $D_{0,n}\times \lbrack 0,1]\subset B_{0}=D^{2}\times \lbrack 0,1],$
for each $n,$

\item $D_{0,n}\times \{0\} \subset D^{2}\times \{0\}$ and $D_{0,n}\times
\{1\} \subset D^{2}\times \{1\},$

\item the 3-balls $H_{1},\cdots ,H_{k}$ are contained in the union $\cup
_{n=1}^{a}(D_{0,n}\times \lbrack 0,1]).$
\end{itemize}

Then, the surface $\widetilde{S}$ is exactly the union of the following
planar surfaces:

\begin{itemize}
\item $S_{0}=D^{2}\times \{0\} \backslash \cup _{n=1}^{a}(D_{0,n}\times
\{0\} \},$

\item $S_{1}=D^{2}\times \{1\} \backslash \cup _{n=1}^{a}(D_{0,n}\times
\{1\} \},$

\item $\cup _{n=1}^{a}S_{0,n},$ where $S_{0,n}=\partial D_{0,n}\times
\lbrack 0,1],$

\item $P=\partial D^{2}\times \lbrack 0,1].$
\end{itemize}

In Figure 10, a disc $D_{0,1}$ and the 3-ball $D_{0,1}\times \lbrack
0,1]\subset D^{2}\times \lbrack 0,1]$ is illustrated in red.

Clearly, cutting $\widetilde{S}$ along the 2-discs $D_{0,n}$ we obtain the
2-sphere $\widetilde{S}_{0}$ which bounds the 3-ball $B_{0}.$

Given this construction, it is straightforward to identify a simple closed
curve $d$ in $\widetilde{S}$ which bounds a disc $D$ in $\widetilde{N}.$
This provides the desired contradiction. In Figure 10, the disc $D$ is
shaded and is illustrated with $d=$ $\partial D\subset (\partial
D_{0,1}\times \lbrack 0,1])\cup (\partial D^{2}\times \lbrack 0,1])\cup
S_{0}\cup S_{1}.$

(2b) $\cup _{i=1}^{m}H_{F_{i}}\subset B_{0}.$

Then we will show that one of the following three cases must hold: either $%
\partial \widetilde{N}$ is compressible in $\widetilde{N},$ or $\widetilde{S}
$ is parallel to $\partial \widetilde{N},$ or $\widetilde{S}$ is
compressible in $\widetilde{N},$ which again contradicts our assumptions.

Indeed, to reconstruct the surface $\widetilde{S},$ we take finitely many
handles $D^{2}\times \lbrack 0,1]$ and attach them to $B_{0},$ by gluing
each $D^{2}\times \{0,1\}$ onto subdiscs $\Delta _{i}$ of $\partial B_{0}=%
\widetilde{S}_{0}.$ This process yields a handlebody $K,$ and it is evident
that $\widetilde{S}=\partial K.$

Similarly, in order to reconstruct $\widetilde{H},$ we attach handles $%
D^{2}\times \lbrack 0,1]$ to the handlebodies $H_{F_{i}},$ gluing each $%
D^{2}\times \{0,1\}$ onto suitable subdiscs $D_{j}$ lying in $\cup
_{i}\partial H_{F_{i}}.$ Moreover, the compression process of $\widetilde{S}$
described above ensures that the handles attached to $\cup _{i}H_{F_{i}}$
lie within the handles attached in $B_{0}.$ It follows that the handlebody $%
\widetilde{H}$ is contained in the handlebody $K,$ i.e. $\widetilde{H}%
\subset K.$

Under these circumstances, we distinguish the following cases:

(i) Some of the components $H_{F_{i}}$ in $B_{0}$ is of genus at least 1.

Then it immediately follows that $\partial \widetilde{H}$ must be
compressible in $\widetilde{N}.$

(ii) All the components $H_{F_{i}}$ in $B_{0}$ are 3-balls and there are
more than one of them.

Then we can identify a simple closed curve $c$ on $\partial B_{0}$ such that 
$c$ bounds a properly embedded 2-disc $D$ in $B_{0}$ which separates the
components $H_{F_{i}}$ and $c$ is an essential curve of $\widetilde{S}.$
Consequently, it follows that $\widetilde{S}$ must be compressible in $%
\widetilde{N}.$

(iii) The family $H_{F_{i}}$ consists of a single component, say $H_{F},$
which is a 3-ball. Then, assuming that $\partial \widetilde{N}$ is not
parallel to $\widetilde{S}$ we can demonstrate through a series of subcases
that $\partial \widetilde{N}$ must be compressible in $\widetilde{N}.$

As a consequence of the preceding analysis, we deduce that $\widetilde{S}%
_{0} $ is incompressible in $\widetilde{M}.$

Now, if the genus of $\widetilde{S}_{0}$ was 1, i.e. $\widetilde{S}_{0}$ is
an incompressible torus then we may deduce the existence of a map $%
T\rightarrow M,$ where $T$ is a torus, which induces a monomorphism at the
level of fundamental groups. This contradicts the assumption that $M$ is
homotopically atoroidal, see the paragraph 2.2 above.

Finally, by projecting $\widetilde{S}_{0}$ in $M$ via the natural projection
map $p:\widetilde{M}\rightarrow M$ we obtain an essential surface $S_{0}$ in 
$M$ which satisfies Haken's conjecture.
\end{proof}

\begin{remark}
\label{connected boundary}If the manifold $\widetilde{N}$ has more that one
components then the filling manifold $\widetilde{M}$ is constructed by
attaching more than one handlebody. In this case, the previous proof of
Theorem \ref{Haken conjecture final} fails in certain places and therefore
cannot be applied to find the incompressible surface $\widetilde{S}_{0}$ in $%
\widetilde{M}.$ For instance, the resulting handlebodies $H_{F_{i}}$ (with
the notation of proof of Theorem \ref{Haken conjecture final})are not
necessarily all contained either in the 3-ball $B_{0}$ or in its complement $%
B_{0}^{c},$ which is crucial to completing the proof. More precisely, the
the existence of curve $b$ in $\widetilde{S}$ that bounds a disc in $%
\widetilde{N}$ is not guaranteed.
\end{remark}

%\begin{a
\noindent{\bf Acknowledgement.}
I wish to express my thanks to Kenichi Oshika for reading a draft of this
article and for all the insightful comments and questions he raised, which
significantly improved its clarity and conciseness. In particular, he
pointed out an error in one of the cases considered in the proof of Theorem %
\ref{attaching handlebody}. I would also like to thank Henry Wilton, who
drew my attention to Corollary \ref{Indeed finite cover}. His crucial
remarks prompted me to include paragraph 6 in this version, thereby
clarifying the preceding corollary.
%\end{acknowledgement}

\bigskip

Charalampos Charitos,

bakis@aua.gr

Agricultural University Athens

Laboratory of mathematics

75, Iera Odos, 11855 Athens, Grece

\end{document}